\documentclass[11pt]{amsart}

\newtheorem{theorem}{Theorem}[section]

\newtheorem{corollary}[theorem]{Corollary}

\theoremstyle{definition}
\newtheorem{definition}[theorem]{Definition}
\newtheorem{example}[theorem]{Example}
\newtheorem{problem}[theorem]{Problem}
\theoremstyle{remark}
\newtheorem{remark}[theorem]{Remark}

\usepackage[margin=1in]{geometry}
\usepackage{amsfonts}
\usepackage{enumerate}
\usepackage{graphicx}
\usepackage{epstopdf}
\usepackage{float}
\usepackage{algorithmic}
\usepackage{hyperref}
\ifpdf
  \DeclareGraphicsExtensions{.eps,.pdf,.png,.jpg}
\else
  \DeclareGraphicsExtensions{.eps}
\fi

\numberwithin{theorem}{section}
\usepackage{caption}
\usepackage{multirow}
\usepackage{amssymb}
\usepackage{cleveref}
\usepackage{url}

\usepackage{amsopn}

\newcommand{\tp }{{\scriptscriptstyle\mathsf{T}}}

\begin{document}

\title{Accurate Solutions of Polynomial Eigenvalue Problems}

\author{Yiling You}
\address{Department of Statistics, 
    University of Chicago, Chicago, IL}
\email{ylyou@uchicago.edu}
\author{Jose Israel Rodriguez}
\address{Computational and Applied Mathematics Initiative, Department of Statistics,
    University of Chicago, Chicago, IL}
\email[corresponding author]{joisro@uchicago.edu}
\author{Lek-Heng Lim}
\address{Computational and Applied Mathematics Initiative, Department of Statistics,
    University of Chicago, Chicago, IL}
\email{lekheng@galton.uchicago.edu}

\begin{abstract}
Quadratic eigenvalue problems (QEP) and more generally polynomial eigenvalue problems (PEP) are among the most common types of nonlinear eigenvalue problems.
Both problems,  especially the QEP, have extensive applications.
A typical approach to solve QEP and PEP is to use a linearization method to reformulate the problem as a higher dimensional linear eigenvalue problem.
In this article, we use homotopy continuation to solve these nonlinear eigenvalue problems without passing to higher dimensions.  
Our main contribution is to show that our method produces substantially more accurate results, and finds all eigenvalues with a certificate of correctness via Smale's $\alpha$-theory. To explain the superior accuracy, we show that the nonlinear eigenvalue problem we solve is better conditioned than its reformulated linear eigenvalue problem, and our homotopy continuation algorithm is more stable than QZ algorithm --- theoretical findings that are borne out by our numerical experiments. Our studies provide yet another illustration of the dictum in numerical analysis that, for reasons of conditioning and stability, it is sometimes better to solve a nonlinear problem directly even when it could be transformed into a linear problem with the same solution mathematically.

\end{abstract}


\maketitle


\section{Introduction}\label{Intro}

The study of polynomial eigenvalue problems (PEPs) is an important topic in numerical linear algebra. 
Such problems arise in partial differential equations and in various  scientific and engineering applications (more on these below).
A special case of PEP that has been extensively and thoroughly studied is
the \textit{quadratic eigenvalue problem} (QEP), often formulated as below,
where we follow the notations in \cite{QEPsurvey}. 
\begin{problem}[QEP]
Let $M,C,K \in \mathbb{C}^{n\times n}$.
The QEP corresponding to these matrices is to determine all solutions $(x,\lambda)$ to the following equation
\begin{equation}\label{QEP_def}
Q(\lambda)x=0\quad \text{ where }\quad Q(\lambda):=\lambda^2 M+\lambda C+K.
\end{equation}
The dimension of the QEP is $n$. We call a solution $(x,\lambda)$ an eigenpair, $\lambda$ an eigenvalue, and $x$ an eigenvector.
\end{problem}
More generally, we may consider  \textit{polynomial eigenvalue problems} (PEP), of which the QEP is the special case when $m = 2$.
\begin{problem}[PEP]\label{problem:pep}
 Let $A_0,A_1,\dots,A_m \in \mathbb{C}^{n\times n}$.
The PEP corresponding to these matrices is to determine all solutions $(x,\lambda)$ to the equation
\begin{equation}\label{eq:pep_def}
P(\lambda)x=0\quad\text{ where }\quad P(\lambda):=\lambda^mA_m+\lambda^{m-1}A_{m-1}+\cdots+A_0.
\end{equation}
Here $m$ is the degree of the PEP and again the solutions are called eigenpairs.
A general PEP has $mn$ eigenpairs (we use `general' in the sense of algebraic geometry; those unfamiliar with this usage may think of it as `random').
\end{problem}
Well-known examples where such problems arise include the acoustic wave problem
 \cite{MIMS_ep2011.116,acoustic_wave}, which gives  a QEP, and the planar waveguide problem
\cite{MIMS_ep2011.116,stowell2010guided}, which gives a PEP of higher degree; they are discussed in \cref{ss:QEPEx} and \cref{ss:PEPEx} respectively.
The encyclopedic survey \cite{QEPsurvey} contains many more examples of QEPs.

There are several methods for solving the QEP, including iterative methods such as Arnoldi method \cite{10.2307/43633863}, Jacobi--Davidson method \cite{Jacobi-Davison1,Jacobi-Davison2}, and linearization method.
However the existing methods invariably suffer from some form of inadequacies: Either (i) they do not apply to PEP of arbitrary degree, or (ii) they require matrices with special structures, or (iii) they only find the largest or smallest eigenvalues. Essentially the only existing method that potentially avoids these inadequacies is the linearization method, as it is based on a reduction to a usual (linear) eigenvalue problem. With this consideration, the linearization method forms the basis for comparison with our proposed method, which does not suffer from any of the aforementioned inadequacies.  We test our method with the software \textsc{Bertini} \cite{bertinibook} and compare  numerical results obtained with those obtained using the linearization method described in \cref{ss:linearization}.

Our main contribution is to propose the use of homotopy method to directly solve a nonlinear problem and find \emph{all} eigenpairs of a PEP or QEP,  suitably adapted to take advantage of the special structures of these problems. Since this article is written primarily for numerical analysts,  we would like to highlight that the study of  \emph{homotopy method for solving a system of multivariate polynomial equations} (as opposed to systems involving non-algebraic or transcendental functions) has undergone enormous progress within the past decade --- both theoretically, with the import of powerful results from complex algebraic geometry, and practically, with the development of new softwares implementing greatly improved algorithms. Some of the recent milestones include the resolution of  Smale's 17th Problem \cite{SIAMNews} and finding all 
tens of millions of solutions of kinematic problems in biologically-inspired linkage design \cite{bioLinkage}. Homotopy method has been used  for symmetric EVPs and GEPs in the traditional, non-algebraic manner; for instance, \cite{Lui-Golub,Zhang-Law-Golub}  rely  on
Raleigh quotient iterations rather than Newton's method  as the corrector method and their results are not certified.

Fortunately for us, PEP and QEP fall in this realm of algebraic problems --- they are special systems of multivariate polynomial equations. In fact, we will see that when formulated as such a system of polynomial equations, the PEP is better conditioned than its alternative formulation as a generalized eigenvalue problem  obtained using linearization. 
Our numerical experiments also show that as the dimension $n$ crosses a threshold of around $20$, homotopy method begins to significantly outperform linearization method in terms of the normwise backward errors. 
 Moreover, our homotopy method approach is not limited to a specific class of matrices but applies robustly to a wide range of matrices ---
 dense 
  sparse, 
  badly scaled, etc. 
In addition, we will see that our homotopy method approach makes it perfect for parallelization. Perhaps most importantly, a unique feature of our approach is that our outputs are certifiable using Smale's $\alpha$-theory. 



\section{Linearization method}\label{ss:linearization}
 A popular approach to solve QEP and PEP is the \textit{linearization method}.
The goal of which is to transform a PEP into a generalized eigenvalue problem (GEP) 
\cite{QEPsurvey} 
involving an equivalent linear $\lambda$-matrix  $A-\lambda B$.
A $2n \times 2n$  linear $\lambda$-matrix is a \textit{linearization} of $Q(\lambda)$ \cite{gohberg2005matrix, lancaster1985theory} if
\begin{equation}
\begin{bmatrix}
Q(\lambda) & 0\\
0 & I_n\\
     \end{bmatrix}
=E(\lambda)(A-\lambda B)F(\lambda)
\end{equation}
where $E(\lambda)$ and $F(\lambda)$ are $2n\times 2n$ $\lambda$-matrices with constant nonzero determinants. The eigenvalues of the quadratic $\lambda$-matrix $Q(\lambda)$ and the linear $\lambda$-matrix $A-\lambda B$ coincide.

Linearizations of PEP are not unique, but two linearizations commonly used in practice are  the \textit{first} and \textit{second companion forms}:
\begin{equation}\label{L1andL2}
\text{L1:} \quad
\begin{bmatrix}
0 & N\\
-K & -C\\
     \end{bmatrix}
-\lambda
\begin{bmatrix}
N & 0\\
0 & M\\
     \end{bmatrix},
\qquad 
\text{L2:} \quad
\begin{bmatrix}
-K & 0\\
0 & N\\
     \end{bmatrix}
-\lambda
\begin{bmatrix}
C & M\\
N & 0\\
     \end{bmatrix},
\end{equation}
where $N$ can be any nonsingular $n\times n$ matrix. The choice between the two companion forms \cref{L1andL2} usually depends on the nonsingularity of $M$ and $K$ \cite{afolabi1987linearization}.



More generally, a PEP~\cref{eq:pep_def} can also be transformed into a GEP 
of dimension $mn$. The most common linearization is called the \textit{companion linearization}
where $A$ and $B$ are:
\begin{equation}\label{eq:gep}
A=
\begin{bmatrix}
A_0 &  &  & & \\
 & I &  &  & \\
 &  & \ddots &  & \\
 &  &  & I & \\
&  & &  & I
\end{bmatrix},\qquad
B=
\begin{bmatrix}
-A_1 & -A_2 & \cdots & \cdots & -A_m\\
I & 0 & \cdots & \cdots & 0\\
0  & I  &  0 & \cdots & 0\\ 
\vdots & \ddots & \ddots & \ddots & \vdots\\
0 & \cdots & 0 & I& 0
\end{bmatrix}.
\end{equation}
We will call such a GEP corresponding to a PEP its \emph{companion GEP}.
One reason for the popularity of the companion linearization is that eigenvectors of PEP can be directly recovered from eigenvectors of this linearization (see  \cite{doi:10.1137/050628283} for details). 

\section{Stability}\label{ss:Accuracy} To assess the stability and quality of the numerical methods in this article, we use the backward error for PEP as defined and discussed in \cite{tisseur2000backward}. Since we do not run into zero or infinite eigenvalues in all the numerical experiments and application problems in this article, our backward errors are always well-defined.
For an approximate eigenpair $(\widetilde{x},\widetilde{\lambda})$ of $Q(\lambda)$, 
the \emph{normwise backward error} is defined as 
\[
\eta(\widetilde{x},\widetilde{\lambda}) := \min \big\{ \epsilon :\bigl(Q(\widetilde{\lambda})+\Delta Q(\widetilde{\lambda})\bigr)\widetilde{x}=0, \;
 \| \Delta M\| \leq \epsilon \alpha_2 ,\; \| \Delta C\| \leq \epsilon \alpha_1,\; \| \Delta K\| \leq \epsilon \alpha_0\bigr\},
\]
where $\Delta Q(\lambda)$ denotes the perturbation 
\[
\Delta Q(\lambda) = \lambda^2 \Delta M+\lambda \Delta C +\Delta K,
\]
and $\alpha_k$'s are nonnegative parameters that allow freedom in how perturbations are measure, e.g., in an absolute sense ($\alpha_k\equiv 1$) or a relative sense ($\alpha_2 = \|M \|$, $\alpha_1 = \|C \|$, $\alpha_0 = \|K \|$).
In \cite{tisseur2000backward}, the backward error is shown to be equal to the following scaled residual,
\begin{equation}\label{scaled residual}
\eta(\widetilde{x},\widetilde{\lambda})=
\frac{\bigl\| Q(\widetilde{\lambda})\widetilde{x}\bigr\|}{\bigl(|\widetilde{\lambda}|^2 \alpha_2+|\widetilde{\lambda}| \alpha_1+\alpha_0\bigr) \|\widetilde{x} \|},
\end{equation}
which has the advantage of being readily computable. 

We adopt usual convention and say that a numerical method is \textit{numerically stable} if all of the computed eigenpairs have backwards errors that are of the same order as the unit roundoff. It is known that a numerically stable reduction of GEP may be  obtained by computing the \textit{generalized Schur decomposition}:
\begin{equation}\label{generalized Schur decomposition}
W^*AZ=S, \qquad W^*BZ=T,
\end{equation}
where $W$ and $Z$ are unitary and $S$ and $T$ are upper triangular. Then $\lambda_i={s_{ii}/t_{ii}}, i=1,\dots, 2n$, with the convention that $s_{ii}/t_{ii}=\infty$ whenever $t_{ii}=0$.

The QZ algorithm \cite{moler1973algorithm,golub2012matrix} is numerically stable for computing the decomposition~\cref{generalized Schur decomposition} and solving the GEP, but it is not stable for the solution of the QEP. To be more precise, one may solve a QEP via the linearization followed by  QZ algorithm applied to resulting GEP. An approximate eigenvector $\widetilde{x}$ to the QEP can be recovered from either the first $n$ components or the last $n$ components of the approximate eigenvector of the GEP computed by the QZ algorithm, a $2n$-vector $\widetilde{\xi}^\tp=(\widetilde{x}^\tp_1,\widetilde{x}^\tp_2)$ that yields the smallest backward error~\cref{scaled residual}. Nevertheless, this method is in general unstable  \cite{tisseur2000backward}.  Even though its backward stability is not guaranteed, it is still the method-of-choice when all  eigenpairs are desired and if the coefficient matrices have no special structure and are not too large.


The discussions above regarding QEP extends to PEP. For an approximate eigenpair $(\widetilde{x},\widetilde{\lambda})$ of $P(\lambda)$,
the normwise backward error is defined in \cite{tisseur2000backward}
as 
\[
\eta(\widetilde{x},\widetilde{\lambda}) := \min \bigl\{ \epsilon :\bigl(P(\widetilde{\lambda})+\Delta P(\widetilde{\lambda})\bigr)\widetilde{x}=0,\;
\| \Delta A_k\| \leq \epsilon \| E_k\|,\; k=0,\dots, m \bigr\}
\]
where $\Delta P(\lambda)$ denotes the perturbation $$\Delta P(\lambda)=\lambda^m \Delta A_m+\lambda^{m-1} \Delta A_{m-1}+\cdots+\Delta A_0$$
and the matrices $E_k, k= 0,\dots,m$ are arbitrary and represent tolerances against which the perturbations $\Delta A_k$ to $A_k$ will be measured. As in the case of QEP, the normwise backward error of PEP may be computed \cite{tisseur2000backward} via the following expression
\begin{equation}\label{PEP scaled residual}
\eta(\widetilde{x},\widetilde{\lambda})=\frac{\bigl\| P(\widetilde{\lambda})\widetilde{x}\bigr\|}{\bigl( \|E_0\| + |\widetilde{\lambda}| \|E_1\| + |\widetilde{\lambda}|^2 \|E_2\| \dots + |\widetilde{\lambda}|^m \|E_m\| \bigr) \|\widetilde{x} \|}.
\end{equation}

\section{Conditioning}\label{ss:Conditioning}

We will see examples in \cref{sec:cond} where a PEP as formulated in \cref{eq:pep_def} is far better conditioned than its companion GEP in \cref{eq:gep} obtained using linearization. 
Thus for accurate solutions, it is better to solve the nonlinear PEP problem directly rather than to first transform it into a mathematically equivalent GEP. As an analogy, this is much like solving a least squares problem directly versus solving its normal equations.

 To characterize the sensitivity of solutions to these problems, we will  describe their condition numbers in this section. As previously mentioned, in all of our numerical experiments and applications, we will not encounter zero or infinite eigenvalues. So in principle we just need a notion of condition number  \cite{tisseur2000backward} that is based on the nonhomogeneous matrix polynomial $P(\lambda)$ in \cref{eq:pep_def}. However, since the \textsc{Matlab} \texttt{polyeig} function that we use for comparison (see \cref{sec:polyeig}) implements a more general notion of condition number \cite{DEDIEU200371, doi:10.1137/050628283}  based on the homogenized version (permitting zero and infinite eigenvalues), we will briefly review the homogenized eigenvalue problem and define its condition number accordingly. 

We rewrite the polynomial matrix $P(\lambda)$ in \cref{eq:pep_def} in homogeneous form
\begin{equation}\label{eq:HPEP}
P(\lambda_0,\lambda_1) = \sum_{i=0}^{m} \lambda_0^i \lambda_1^{m-i} A_i
\end{equation}
and consider eigenvalues as pairs $(\lambda_0,\lambda_1) \neq (0,0)$ that are solutions of the equation $\det P(\lambda_0,\lambda_1)=0$. Let $T_{(\lambda_0,\lambda_1)}\mathbb{P}_1$ denote the tangent space at $(\lambda_0,\lambda_1)$ to $\mathbb{P}_1$, the projective space of lines through the origin in $\mathbb{C}^2$. A condition operator is defined in \cite{DEDIEU200371} as $K(\lambda_0,\lambda_1):(\mathbb{C}^{n \times n})^{m+1} \to T_{(\lambda_0,\lambda_1)}\mathbb{P}_1$ of the eigenvalue $(\lambda_0,\lambda_1)$ as the differential of the map from the $(m+1)$-tuple $(A_0,\dots,A_m)$ to $(\lambda_0,\lambda_1)$ in projective space. If we write a representative of an eigenvalue $(\lambda_0,\lambda_1)$ as a row vector $[\lambda_0,\lambda_1] \in \mathbb{C}^{1 \times 2}$, the condition number $\kappa_P (\lambda_0,\lambda_1)$ can be defined as a norm of the condition operator \cite{doi:10.1137/050628283},
\[
\kappa_P(\lambda_0,\lambda_1) := \max_{\| \Delta A\|\leq 1}\frac{\|K(\lambda_0,\lambda_1)\Delta A\|_2}{\|[\lambda_0,\lambda_1]\|_2},
\]
for any arbitrary norm on $\Delta A$. We will choose the norm on $(\mathbb{C}^{n \times n})^{m+1}$ to be the $\omega$-weighted Frobenius norm
\[
\|A\| :=\|(A_0,\dots,A_m)\|=\bigl\|[\omega_0^{-1}A_0,\dots,\omega_m^{-1}A_m]\bigr\|_F,
\]
with weights $\omega_i > 0$, $i =1,\dots,m$. If we define the operators $\partial_{\lambda_0}:= \partial/\partial \lambda_0$ and $\partial_{\lambda_1}:= \partial/\partial \lambda_1$, then the normwise condition number $\kappa_P(\lambda_0,\lambda_1)$ of a simple eigenvalue $(\lambda_0,\lambda_1)$ is given by
\[
\kappa_P(\lambda_0,\lambda_1) = \left( \sum_{i=0}^m |\lambda_0|^{2i}|\lambda_1|^{2(m-i)}\omega_i^2 \right)^{1/2} \frac{\|y\|_2\|x\|_2}{|y^*(\overline{\lambda}_1\partial_{\lambda_0}P-\overline{\lambda}_0\partial_{\lambda_1}P)|_{(\lambda_0,\lambda_1)}x|},
\]
where $x,y$ are the corresponding right and left eigenvector respectively \cite{doi:10.1137/050628283}.








\begin{remark}
Homogenization allows one to better handle eigenvalues at infinity.
In general, the characteristic polynomial is the determinant of the matrix polynomial and takes the form  $\det (A_m)\lambda_1^{mn}+\cdots+\det (A_0)\lambda_0^{mn}$.
Therefore, $P(\lambda_0,\lambda_1)$ has $mn$ finite eigenvalues when the matrix $A_m$ is nonsingular.
However, when $\det(A_m)=0$, the characteristic polynomial has degree $r < mn$ and there are $r$ finite eigenvalues and $mn-r$ infinite eigenvalues.
Those infinite eigenvalues correspond to $\lambda_0=0$. None of the numerical experiments we consider in this article has eigenvalues at infinity.
\end{remark}

\section{Certification}\label{sec:cert}
 
A major advantage of our homotopy method approach for solving the PEP is that we may use Smale's $\alpha$-theory (also known as Shub--Smale's $\alpha$-theory, see \cite[Chapter 8]{BCSS1998}) to certify that the Newton iterations will converge quadratically to an eigenpair. In numerical analysis lingo, this means we can control the \emph{forward error}, not just the backward error.

To apply Smale's $\alpha$-theory, we view the PEP as a collection of $n$ polynomial equations  \cref{eq:pep_def} 
and an affine linear constraint $L(x)=0$, which yields a polynomial map $f = (f_0,\dots,f_n):\mathbb{C}^{n+1}\to\mathbb{C}^{n+1}$.
The affine linear polynomial $L(x)$ is chosen randomly as described in the numerical experiments 
\cref{Numerical results}.


Let $\mathcal{V}(f):=\left\{
\zeta\in\mathbb{C}^{n+1} : f(\zeta)=0
\right\}$ and let $Df(z)$ be the Jacobian matrix of the system $f$ at $z=(x,\lambda)$.
Consider the map $N_f:\mathbb{C}^{n+1}\to\mathbb{C}^{n+1}$ defined by 
\[
N_f(z):=\begin{cases}
z-Df\left(z\right)^{-1}f\left(z\right) & \text{if }Df\left(z\right)\text{ is invertible,}\\
z & \text{otherwise}.
\end{cases}
\]
We say the point $N_f(z)$ is the \textit{Newton iteration of $f$ starting at $z$}. 
The $k$th Newton iteration of $f$ starting at $x$ is denoted by $N_f^{k}(z)$. 
Now we define precisely what we mean by an approximate solution to $f$. 

\begin{definition}\label{def:approxSol}\emph{\cite[p.~155]{BCSS1998}}
With the notation above, a point $z\in\mathbb{C}^{n+1}$ is an \textit{approximate solution to $f$ with \textit{associated solution} $\zeta\in\mathcal{V}(f)$}, if for every $k\in\mathbb{N}$, 
\begin{equation}
 \left\Vert N_{f}^{k}\left(z\right)-\zeta\right\Vert \leq\left(\frac{1}{2}\right)^{2^{k}-1}\left\Vert z-\zeta\right\Vert,
\end{equation}
where the norm is the $2$-norm
$\lVert z
\rVert=\bigl(
|z_1|^2+\cdots+|z_{n+1}|^2
\bigr)^{1/2}$.
\end{definition}

Smale's $\alpha$-theory gives a condition for when a given point $z$ is an approximate solution to $f=0$ using the following constants when $Df(z)$ is invertible:
\begin{align*}
\alpha(f,z) &:=\beta(f,z)\gamma(f,z),\\
\beta(f,z)  &:=\left\Vert z-N_f(z) \right\Vert=
\left\Vert
Df(z)^{-1}f(z)
\right\Vert,\\
\gamma(f,z) &:=
\sup_{k\geq 2}\left\Vert
\frac{Df(z)^{-1}D^kf(z)}{k!}
\right\Vert^{1/(k-1)}.
\end{align*}
The following theorem from \cite{HS12} is a version of Theorem~2 in \cite[p.~160]{BCSS1998} and it  provides a certificate  
that a point $z$ is an approximate solution to $f=0$.

\begin{theorem}\label{smale}
If $f:\mathbb{C}^{n+1}\to\mathbb{C}^{n+1}$ is a polynomial system and $z\in\mathbb{C}^{n+1}$,
with 
\[
\alpha(f,z)<\frac{13-3\sqrt{17}}{4}\approx 0.157671,
\]
then $z$ is an approximate solution to $f=0$.
\end{theorem}

The quantity $\gamma(f,z)$ is difficult to compute in general. 
In \cite{HS12},  this quantity is bounded in terms of  an alternative that is more readily computable.
Define 
\begin{equation}\label{eq:mu}
\mu(f,z):=
\max\left\{
1,\Vert f\Vert\cdot \Vert Df(z)^{-1}
\Delta_{(d)}(z)\Vert
\right\},
\end{equation}
where
$\Delta_{(d)}(z)$ 
is an $(n+1)\times(n+1)$ diagonal matrix with 
$i$th diagonal entry $d_i^{1/2}(1+\Vert z\Vert^2)^{(d_i-1)/2}$ and $d_i := \operatorname{deg}(f_i)$. 
For a polynomial $g=\sum_{|\nu|\leq d}a_\nu z^\nu$, we define the norm $\Vert g\Vert$ according to \cite{HS12},
\begin{equation}\label{eq:norm1}
\Vert g\Vert^2:=
\sum_{|\nu|\leq d}  |a_\nu|^2\frac{\nu!(d-|\nu|)!}{d!}.
\end{equation}
This can extended to polynomial system $f$ via
\begin{equation}\label{eq:norm2}
\Vert f\Vert^2:=\sum_{i=0}^{n}\Vert f_i\Vert^2
\end{equation}
where each $ \Vert f_i\Vert^2$ is as defined in \cref{eq:norm1}.

With the notations above, \cite{HS12} gives the following bound:
\begin{equation}\label{eq:gammaBound}
\gamma(f,z)\leq \frac{\mu(f,z)d_{\max}^{3/2}}{2(1+\Vert z\Vert^2)^{1/2}},
\end{equation}
where $d_{\max}$ is the maximal degree of a polynomial in the system. 

When  $f$ comes from a PEP,  \cref{pep_certification} follows from \cref{eq:mu}--\cref{eq:norm2} and a straightforward
calculation of $Df(z)^{-1}\Delta_{(d)}(z)$. This yields the value of $\mu(f,z)$ and thus a bound for $\gamma(f,z)$ for the PEP.
We will rely on our bound in the acoustic wave problem \cref{Certification_APP1} to certify the eigenvalues of the PEP therein. 

\begin{corollary}\label{pep_certification}
Let $f=0$ denote the PEP in \cref{problem:pep}, and 
let $L$ denote the affine constraint on the $x$.  
Then,
\begin{equation}\label{mu_formula}
\mu(f,x,\lambda) =
\left[\Vert L\Vert^2+\sum_{k=0}^m \frac{k!(m-k)!}{(m+1)!}\|A_k\|_F^2\right]^{1/2} 
\left\|
\begin{bmatrix}
\frac{P(\lambda)}{\sqrt{m+1}(1+\|[x^\tp,\lambda]\|^2)^{m/2}} & P'(\lambda)x  \\
\frac{(\nabla_x L)^\tp}{\sqrt{m+1}(1+\|[x^\tp,\lambda]\|^2)^{m/2}} & 0 
\end{bmatrix}
^{-1}\right\|.
\end{equation}
\end{corollary}


\section{Homotopy method for the polynomial eigenvalue problem}\label{ss:homotopymethod}

We now describe our approach of using  \emph{homotopy continuation method}, also called \emph{homotopy method}, to solve PEPs.
Our goal is to find all eigenpairs.

Homotopy method deforms solutions of a
\emph{start system} $S(z)=0$ that is easy to solve 
to solutions of a 
\emph{target system} $T(z)=0$ 
that is of interest. 
More precisely, 
a \emph{straight-line homotopy}
with \emph{path parameter} $t$ is defined as 
\begin{equation}\label{eq:H}
H(z,t):=(1-t)S(z)+tT(z), \quad t\in[0,1].
\end{equation}
When $t=0$ or $t=1$,
the system $H(z,t) = 0$ gives the start system
 $H(z,0)=S(z)=0$ or the target system $H(z,1)=T(z)=0$ respectively.

A start system for the homotopy $\eqref{eq:H}$
is said to be \emph{chosen correctly} if the 
 following properties \cite{li_1997} hold:

\begin{enumerate}[\upshape (i)]
\item the solution set of the start system 
$S(z)=0$ are known or easy to obtain;
\item the solution set of 
$H(z,t)=0$ for $0 \leq t < 1$ consists of a finite number of smooth paths, each parametrized by $t$ in $[0,1)$;
\item for each isolated solution of the target systems
$T(z)=0$, there is some path originating at $t=0$, that is, a solution of the start systems 
$S(z)=0$.
\end{enumerate}

\begin{example}
For an example of  homotopy familiar to numerical linear algebraists, consider $P(\lambda) = I - \lambda A$, i.e., $m=1$, $A_0 = I$, and $A_1 = -A$. 
Let $D$ be the diagonal matrix whose diagonal entries are the diagonal entries of $A$. 
The proof of strengthened Gershgorin Circle Theorem \cite{varga} exactly illustrates a   straight-line homotopy path $H(t)=(1-t)D+tA$, $t\in[0,1]$. Note however that such a $D$ would be a poor choice for us as the start system for the homotopy is not guaranteed to be chosen correctly: The solution set  $H(t)=0$ for $0 \leq t < 1$ need not consist of a finite number of smooth paths. 
\end{example}

In this article, we consider a target system to solve 
the PEP in \eqref{eq:pep_def} given by
\begin{equation}\label{eq:targetpep}
T(z)=\begin{bmatrix}
P(\lambda)x\\
L(x)
\end{bmatrix}
\end{equation}
 where  $z=(x,\lambda)$ and
  $L(x)$  is a general affine linear polynomial, 
chosen randomly so that we have a polynomial system $T : \mathbb{C}^{n+1} \to \mathbb{C}^{n+1}$ as defined in \cref{sec:cert}. The requirement that $L(x) = 0$ also fixes the scaling indeterminacy in the polynomial eigenvector $x$.  In more geometric language, $x$ is a point in projective space and  the random choice of  $L(x)$ restricts this space to a general affine chart 
so that the eigenvectors are not at infinity.  If one instead had chosen $L(x)=x_1-1$, then eigenvectors with a first coordinate equal to zero would not be solutions to the~system.

There is an obvious choice of start system
--- 
we choose random diagonal matrices $D_i$'s to replace  the coefficient matrices $A_i$'s in $P(\lambda)$:
\begin{equation}\label{eq:startpep}
S(z)=\begin{bmatrix}
(\lambda^mD_m+\lambda^{m-1}D_{m-1}+\cdots+
D_0)x\\
L(x)
\end{bmatrix}
\end{equation}

One observes that  $S(x,\lambda)=0$ is a polynomial system with linear products of $(x,\lambda)$.
Specifically, let $d_{j,i}$ denote the $i$th diagonal entry of $D_j$ and let $x_i$ be the $i$th entry of $x$. Then we may factor the univariate polynomial
\[
\bigl(\lambda^m d_{m,i}+\lambda^{m-1} d_{m-1,i}+\dots+d_{0,i}\bigr)x_i=d_{m,i}(\lambda-r_{m,i})( \lambda-r_{m-1,i}) \cdots ( \lambda-r_{1,i})x_i,
\]
where $ r_{j,i}\in \mathbb{C}$, $j =1,\dots, m$, are the roots of the respective monic polynomial (obtained via, say, the Schur decomposition of its companion matrix).
The solutions to \eqref{eq:startpep} are then simply
\[
\lambda = r_{j,i}, \qquad
x_k = 0 \; \text{for all}\; k \neq i, \qquad
L(0,\dots,0,x_i,0,\dots,0)=0,
\]
for every $j = 1,2,\dots,m$ and $i=1,2,\dots,n$.

From the known eigenpairs at $t=t_0$, solutions at $t=t_0+\Delta t$ can be obtained by iterative methods \cite{allgower2012numerical,allgower1993continuation}. This process is called path tracking. The steps in iteration  are repeated until $t$ reaches $1$ or $t$ is sufficiently close to $1$ by some criteria. The output is regarded as an approximation of a solution to the PEP.

The path tracking becomes numerically unstable as we get close to the space of ill-posed systems,
also called the \textit{discriminant locus}. 
In other words, this occurs when the Jacobian of $H(z,t_0)$ with respect to $z$ has large condition number. 
Using randomization and adaptive-precision we  are able to avoid these ill-conditioned problems in the middle of path tracking.
This is because the discriminant locus has real codimension two in the space of coefficients, whereas the homotopy is over a real one-dimensional path.

\begin{example}
A QEP may be expressed as a univariate polynomial root finding problem by taking the determinant of $Q(\lambda)$. 
Say, $n=2$, then we have 
\[
f(\lambda,M,C,K):=\det Q(\lambda)=\det\biggl(\lambda^2
\begin{bmatrix}
      m_{11} &  m_{12}\\
      m_{21} & m_{22}\\
     \end{bmatrix}
+\lambda
\begin{bmatrix}
      c_{11} &  c_{12}\\
      c_{21} & c_{22}\\
     \end{bmatrix}
+
\begin{bmatrix}
      k_{11} &  k_{12}\\
      k_{21} & k_{22}\\
     \end{bmatrix}
\biggr).
\]
Solving this univariate problem is numerically unstable --- in fact this is a special case of the resultant method, which is unstable in general \cite{NT2016}.
The discriminant locus is the set of $(M,C,K) \in \mathbb{C}^{2 \times 2 \times 3}$ is such that there exists a $\lambda$ where the following is satisfied:
 $f(\lambda,M,C,K)= {\partial f}/{\partial \lambda}=0.$

If we were to express $f$ as a sum of monomials then it is the determinant of a Sylvester matrix. Recall that the entries of a Sylvester matrix of two polynomials are coefficients of the polynomials, which in our case are $f$ and $\partial f/\partial \lambda$, i.e.,
\[
{
\begin{bmatrix}
      s_4 & s_3 &s_2&s_1&s_0& &\\
      &s_4 & s_3 &s_2&s_1&s_0& \\
& &s_4 & s_3 &s_2&s_1&s_0\\
4s_4&3s_3&2s_2&s_1& & &\\
&4s_4&3s_3&2s_2&s_1& &\\
&&4s_4&3s_3&2s_2&s_1& \\
&&&4s_4&3s_3&2s_2&s_1\\
     \end{bmatrix},
}
\]
where $s_i$'s are polynomials in the entries of $(M,C,K)$.
The determinant of this matrix is zero when the two polynomials have a common root, i.e.,
\[
f(\lambda,M,C,K)=s_4\cdot\lambda^4+s_3\cdot\lambda^3+s_2\cdot\lambda^2+s_1\cdot\lambda+s_0.
\]
So the vanishing of this determinant defines where the problem is ill-conditioned. By taking a randomized start system, the homotopy method will avoid this space.
\end{example}

\section{Numerical experiments}\label{Numerical results}

In this section, we provide numerical experiments to compare the speed and accuracy of (i) solving the PEP as formulated in \cref{eq:pep_def} with homotopy continuation (henceforth abbreviated as \emph{homotopy method}), and (ii) solving the companion GEP  in \cref{eq:gep} with QZ algorithm (henceforth abbreviated as \emph{linearization method}). We will also compare the conditioning of the PEP in \cref{eq:pep_def} and its companion GEP in \cref{eq:gep}.  All experiments were performed on a computer running Windows 8 with an Intel Core i7 processor. For the linearization method, we use the implementation in the \textsc{Matlab} \texttt{polyeig} function; for the homotopy method, we implemented it with \textsc{Bertini} 1.5.1. We tested both methods  in serial on a common platform: \textsc{Matlab} R2016a with the interface \textsc{BertiniLab} 1.5 \cite{Bates2016}. We also tested the homotopy method in parallel using Intel MPI Library 5.1 compiled with intel/16.0 on the Midway1 compute cluster\footnote{\url{https://rcc.uchicago.edu/resources/high-performance-computing}} in the University of Chicago Research Computing Center.

\subsection{The \texttt{polyeig} and  \texttt{quadeig} functions}\label{sec:polyeig} In our numerical experiments,  we rely on the \texttt{polyeig} function in  \textsc{Matlab} to solve PEPs with linearization method and compute condition numbers. This routine adopts the companion linearization, uses QZ factorization to compute generalized Schur decomposition, recovers the right eigenvectors that has the minimal residual,  and gives the condition number for each eigenvalue.  While different linearization forms can have widely varying eigenvalue condition numbers \cite{doi:10.1137/050628283}, \texttt{polyeig} does not allow for other choices.

A more recent linearization-based algorithm, \texttt{quadeig}, was proposed in \cite{Hammarling:2013:ACS:2450153.2450156} for QEPs. While \texttt{polyeig} is a built-in \textsc{Matlab} function, \texttt{quadeig} is a third-party program implemented in \textsc{Matlab}. It incorporates extra preprocessing steps that scale the problem's  parameters and choose the linearization with favorable conditioning and backward stability properties. However, we remain in favor of using \texttt{polyeig} as the basis of comparison against our homotopy method for solving QEPs as
(i) the scaling is redundant since the $2$-norms of our random coefficient matrices are all approximately one;
(ii) for a uniform comparison, we prefer a method that works for all $m$ but \texttt{quadeig} does not extend to higher order PEP since the scaling process cannot be generalized. There will be further discussions of these issues in \cref{sec:cond}.

%


\subsection{Speed comparisons for QEP}\label{sec:speedQEP} We first test and compare timings for our homotopy method and the linearization method for general quadratic matrix polynomials. The matrix polynomial is generated at random\footnote{Results on matrix polynomials arising from actual applications will be presented in \cref{sec:app}.}  with coefficient matrices $M,C,K$ in \cref{eq:pep_def} having independent entries following the standard complex Gaussian distribution $\mathcal{N}_{\mathbb{C}}(0,1)$, using the \texttt{randn} function in \textsc{Matlab}. The coefficients of $L(x)$ are also chosen randomly in the same way. We performed three sets of numerical experiments:
\begin{enumerate}[\upshape (i)]
\item\label{item:lin}
\cref{timing_linearization} gives the elapsed timings for computing eigenpairs with the linearization method for dimensions $n=2,\dots,100$.
\item\label{item:serial}
\cref{timing_homotopy_serial} gives the elapsed timings for computing eigenpairs with the homotopy method for dimensions $n=2,\dots,80$.
\item\label{item:par}
\cref{timing_homotopy_parallel} gives the elapsed timings for computing the eigenpairs with homotopy method  in parallel on $20$ cores for dimensions $n=20, 30\dots,100$. 
\end{enumerate}
For \cref{item:lin} and \cref{item:serial} the timings include both the setup and the solution process; for \cref{item:par}, the timings only include the solution process running \textsc{Bertini} in parallel. For each method and each dimension, we run our experiments ten times and record the best, average, median, and worst performance. The conclusions drawn from these experiments are discussed in \cref{sec:conclude}.

\begin{table}[H]
\footnotesize
\centering
\caption{\textsc{Speed of QEP --- linearization.} Elapsed timings (in seconds) for the linearization method with dimensions $n=2,3,\dots,100$ (not all displayed).}
\label{timing_linearization}
\begin{tabular}{ c| c| c| c| c|c}
  $n$ & \# ROOTS & BEST & AVERAGE & MEDIAN & WORST \\
\hline
2 & 4 & 2.4550E-4&	0.0015	&8.0279E-4&	0.0082 \\
3 & 6 & 3.1406E-4&	6.7019E-4	&6.0923E-4&	0.0013 \\
4 & 8 & 3.6414E-4&	6.7064E-4	&5.9055E-4&	0.0011 \\
5 & 10 & 4.2818E-4&	6.7003E-4	&5.2938E-4&	0.0013 \\
6 & 12 & 4.6144E-4&	0.0010	&7.3403E-4&	0.0034 \\
7 & 14 & 5.8542E-4&	9.3642E-4	&9.1938E-4&	0.0015 \\
8 & 16 & 6.7327E-4&	9.5284E-4	&8.5719E-4&	0.0013 \\
9 & 18 & 7.3280E-4&	0.0027	&0.0012&	0.0158 \\
10 & 20 & 9.2533E-4&	0.0015	&0.0013&	0.0028 \\
11 & 22 & 9.2246E-4&	0.0014	&0.0012&	0.0027 \\
12 & 24 & 0.0011&	0.0017	&0.0015&	0.0027 \\
13 & 26 & 0.0012&	0.0017	&0.0017&	0.0025 \\
14 & 28 & 0.0016&	0.0020	&0.0019&	0.0025 \\
15 & 30 & 0.0015&	0.0023	&0.0021&	0.0035 \\
16 & 32 & 0.0020&	0.0030	&0.0029&	0.0050 \\
17 & 34 & 0.0025&	0.0037	&0.0033&	0.0058 \\
18 & 36 & 0.0029&	0.0037	&0.0036&	0.0056 \\
19 & 38 & 0.0032&	0.0048	&0.0040&	0.0091 \\
\hline
20 & 40 & 0.0031&	0.0038	&0.0038&	0.0048 \\
30 & 60 & 0.0089&	0.0109	&0.0103&	0.0164 \\
40 & 80 & 0.0114&	0.0157	&0.0152&	0.0220 \\
50 & 100 & 0.0225&	0.0253	&0.0252&	0.0295 \\
60 & 120 & 0.0328&	0.0407	&0.0400&	0.0511 \\
70 & 140 & 0.0418&	0.0527	&0.0538&	0.0624 \\
80 & 160 & 0.0642&	0.0762	&0.0768&	0.0850 \\
90 & 180 & 0.0867&	0.0983	&0.0957&	0.1191 \\
100 & 200 & 0.1177&	0.1285	&0.1229&	0.1504
\end{tabular}
\end{table}
\begin{table}[H]
\footnotesize
\centering
\caption{\textsc{Speed of QEP --- homotopy in serial.}
Elapsed timings (in seconds) for the homotopy method with dimensions $n=2,3,\dots,80$ (not all displayed).}
\label{timing_homotopy_serial}
\begin{tabular}{ c| c| c| c| c|c}
  $n$ & \# ROOTS &  BEST & AVERAGE & MEDIAN & WORST \\
\hline
2 & 4 & 0.3692&	0.4185	&0.3950&	0.6064 \\
3 & 6 & 0.4717&	0.5509	&0.5225&	0.7905 \\
4 & 8 & 0.5005&	1.1083	&0.5553&	5.9034 \\
5 & 10 & 0.6316&	0.7809	&0.6735&	1.0826 \\
6 & 12 & 0.7767&	1.0083	&0.9946&	1.2599 \\
7 & 14 & 0.8593&	1.0333	&0.9418&	1.4684  \\
8 & 16 & 1.0043&	1.4300	&1.1941&	2.6119  \\
9 & 18 & 1.2054&	1.7269	&1.3505&	2.8268  \\
10 & 20 & 1.4373&	1.6836	&1.5428&	2.9608  \\
11 & 22 & 1.6972&	2.2203	&2.0245&	3.3329  \\
12 & 24 & 1.9914&	2.7784	&2.3702&	4.6917  \\
13 & 26 & 2.1979&	3.3690	&2.6251&	7.5556  \\
14 & 28 & 2.6687&	3.8757	&3.7469&	5.3345  \\
15 & 30 & 3.1057&	5.2769	&5.1792&	9.3289  \\
16 & 32 & 3.4994&	7.2866	&6.6201&	16.3712  \\
17 & 34 & 3.8719&	6.1690	&5.9380&	9.3421  \\
18 & 36 & 4.9051&	7.8015	&8.2610&	13.4794  \\
19 & 38 & 5.4613&	10.4376	&11.2833&	13.4682  \\
\hline
20 & 40 & 6.4268&	14.0531	&14.0381&	23.9232  \\
30 & 60 & 15.2266&	32.1909	&32.1480&	63.8345  \\
40 & 80 & 74.4945&	103.9896	&93.9979&	161.6461  \\
50 & 100 & 133.0395&	244.3541	&245.9194&	394.0840  \\
60 & 120 & 309.0921&	532.1712	&485.0595&	943.0120  \\
70 & 140 & 705.9720&	1200.4053	&1101.9838&	1796.4483  \\
80 & 160 & 1207.1300&	1848.9342	&1745.6093&	2974.4295
\end{tabular}
\end{table}
\begin{table}[H]
\footnotesize
\centering
\caption{\textsc{Speed of QEP --- homotopy in parallel.}
Elapsed timings (in seconds) for the homotopy method ran in parallel on $20$ cores and with dimensions $n=20, 30,\dots,100$.}
\label{timing_homotopy_parallel}
\begin{tabular}{ c| c| c| c| c|c}
  $n$ & \# ROOTS &  BEST & AVERAGE & MEDIAN & WORST \\
\hline
20 & 40 & 2.5010&	2.8091	&2.7410&	3.5160 \\
30 & 60 & 2.5890&	5.6856	&4.6965&	9.4390 \\
40 & 80 & 8.3950&	12.7100	&13.2490&	16.3700 \\
50 & 100 & 16.1620&	22.2175	&21.8220&	26.9460 \\
60 & 120 & 18.2790&	24.4677	&23.4870&	29.1530 \\
70 & 140 & 36.5210&	53.7870	&54.7050&	64.7610  \\
80 & 160 & 63.6320&	83.7379	&77.9380&	105.2450  \\
90 & 180 & 101.4390&	142.3974	&149.6170&	177.7420  \\
100 & 200 & 164.3740&	196.5620	&199.4585&	239.5180
\end{tabular}
\end{table}

\subsection{Accuracy comparisons for QEP}\label{sec:accuracyQEP} 
\begin{figure}[H]
\centering
\caption{Graphs for \cref{Ave_abs_berr}. Blue for linearization and red for homotopy.}
\includegraphics[scale=0.36]{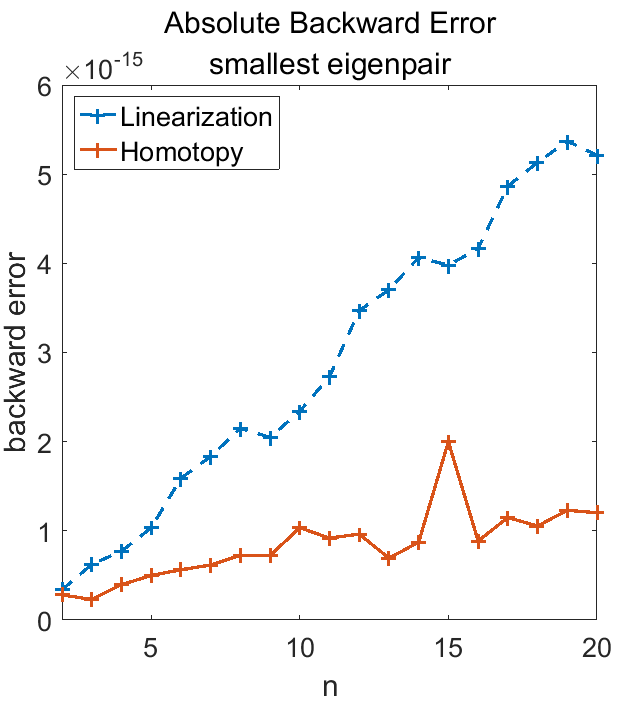}
\includegraphics[scale=0.36]{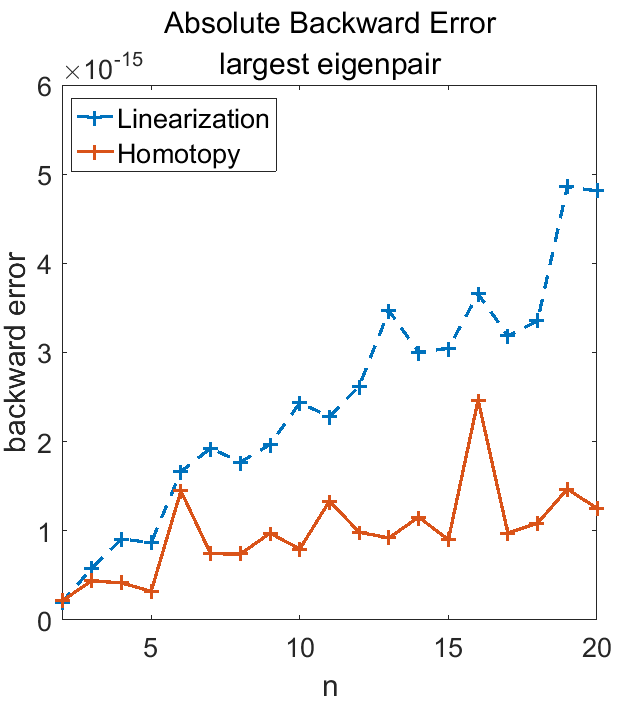}\\
\includegraphics[scale=0.36]{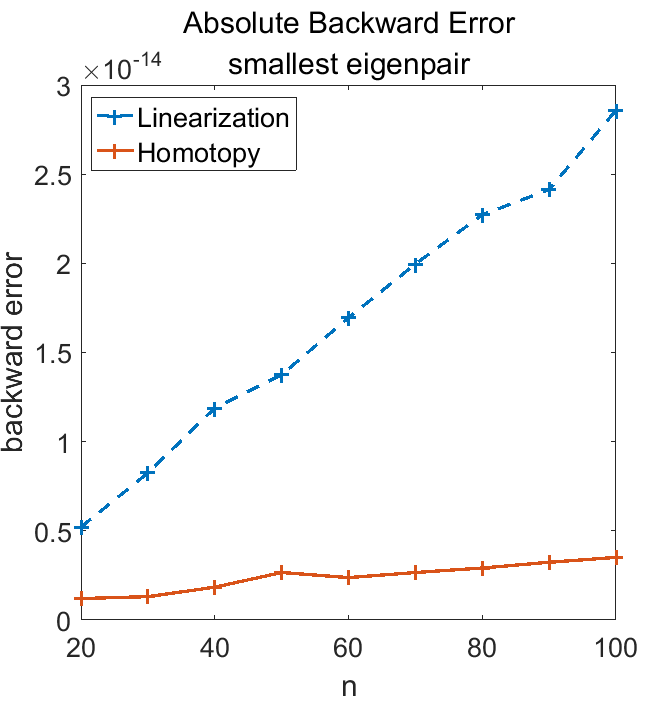}
\includegraphics[scale=0.36]{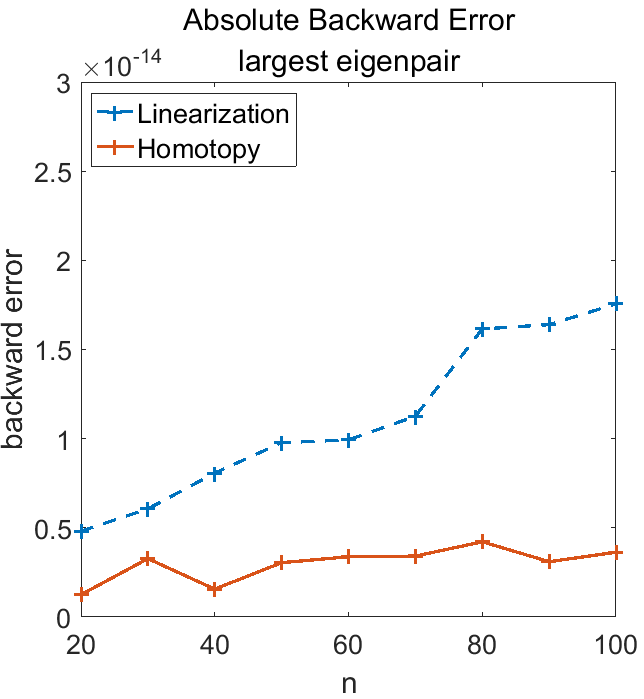}
\label{Ab_BErr_2-100}
\end{figure}

\begin{figure}[H]
\caption{Graphs for \cref{Ave_rel_berr}.
Blue for linearization and red for homotopy.}
\centering
\includegraphics[scale=0.36]{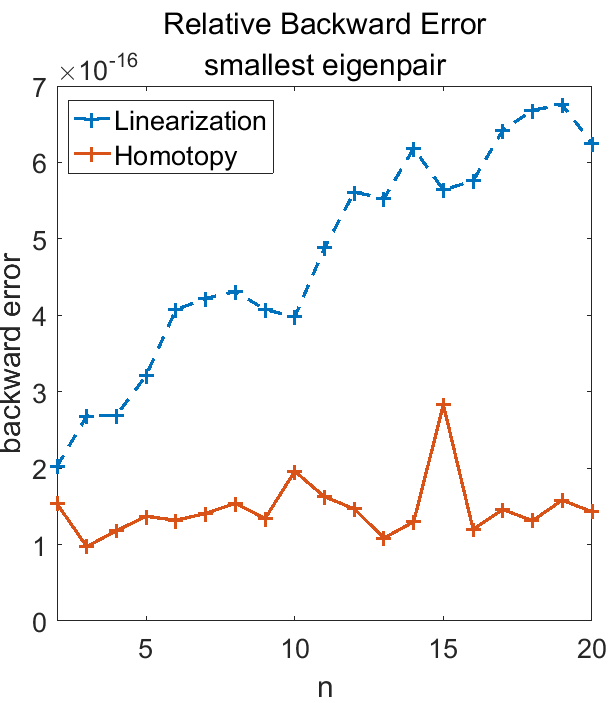}
\includegraphics[scale=0.36]{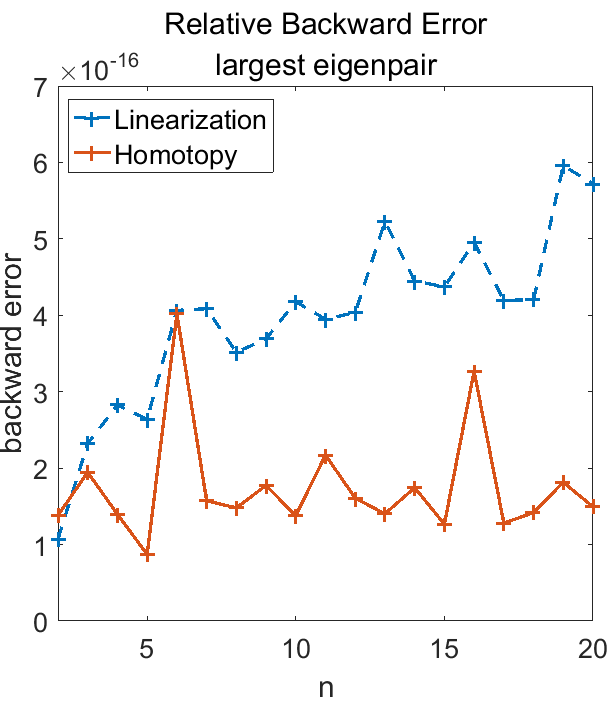}\\
\includegraphics[scale=0.36]{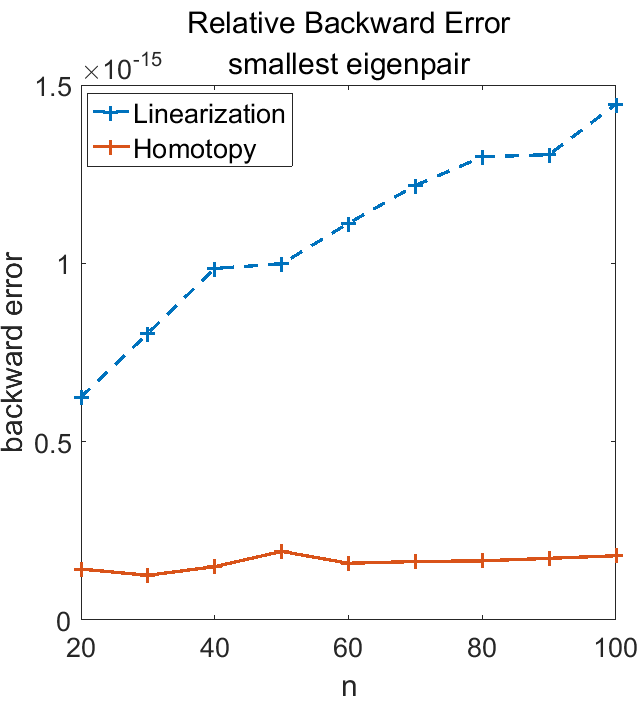}
\includegraphics[scale=0.36]{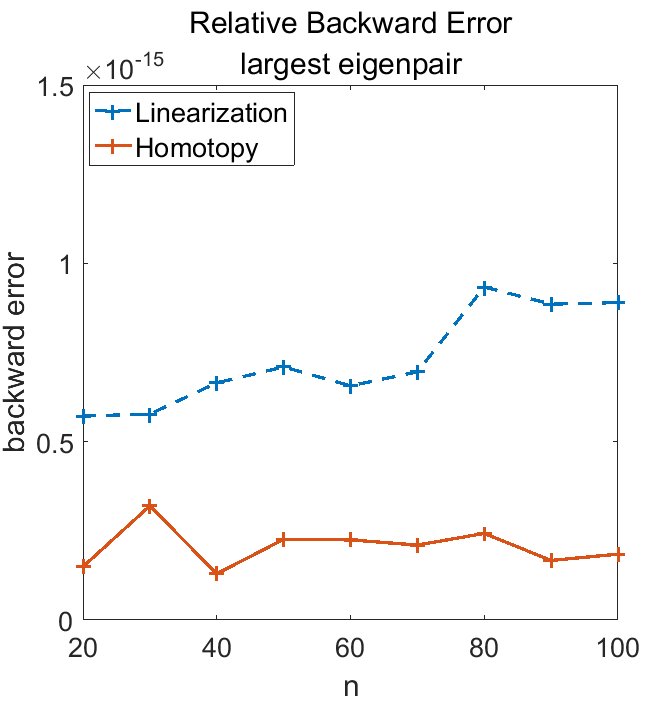}
\label{Rel_BErr_2-100}
\end{figure}

Next we test and compare the absolute and relative backward errors of the computed eigenpairs corresponding to the smallest and largest eigenvalues for our homotopy method and the linearization method for randomly generated matrix polynomials.  For each of the methods, we tabulated the absolute and relative backward errors of the computed eigenpairs corresponding to the smallest and largest eigenvalues  with dimension $n=2,\dots,100$. 
All tests are averaged over ten runs. They are compared side-by-side in \cref{Ave_abs_berr} (absolute error) and \cref{Ave_rel_berr} (relative error), and graphically in \cref{Ab_BErr_2-100} (absolute error) and \cref{Rel_BErr_2-100} (relative error).
The conclusions drawn from these experiments are presented in \cref{sec:conclude}.

\begin{table}[H]
\footnotesize
\centering
\caption{\textsc{Accuracy of QEP --- linearization vs homotopy.}
Absolute backward errors of computed smallest and largest eigenpairs with dimensions $n=2,\dots,100$. }
\label{Ave_abs_berr}
\begin{tabular}{ c| c| c| c| c}
 &\multicolumn{2}{r|}{LINEARIZATION AVG.\ BK.\ ERR.}& \multicolumn{2}{c}{HOMOTOPY AVG.\ BK.\ ERR.} \\
\cline{2-5}
  $n$ &SMALLEST&LARGEST& SMALLEST&LARGEST  \\
\hline
2 & 3.39024E-16&2.01468E-16 & 2.81624E-16&2.13291E-16 \\
5 & 1.03302E-15&8.68237E-16 & 4.98414E-16&3.21280E-16	 \\
10 & 2.34081E-15&2.43521E-15 & 1.03488E-15&7.96172E-16	 \\
20 & 5.21350E-15&4.81234E-15 & 1.20237E-15&1.25302E-15	 \\
30 & 8.24479E-15&6.06398E-15 & 1.31006E-15&3.27204E-15	 \\
40 & 1.18830E-14&8.06968E-15 & 1.83843E-15&1.56400E-15	 \\
50 & 1.37506E-14&9.77673E-15 & 2.66411E-15&3.04532E-15	 \\
60 & 1.69521E-14&9.92973E-15 & 2.37774E-15&3.37814E-15	 \\
70 & 1.99600E-14&1.12571E-14 & 2.66186E-15&3.414932E15	 \\
80 & 2.27208E-14&1.61510E-14 & 2.91539E-15&4.21745E-15	 \\
90 & 2.41366E-14&1.63954E-14 &	3.23766E-15&3.10200E-15 \\
100 & 2.85801E-14&1.75552E-14 &	3.51036E-15&3.62939E-15
\end{tabular}
\end{table}

\begin{table}[H]
\footnotesize
\centering
\caption{\textsc{Accuracy of QEP --- linearization vs homotopy.}
Relative backward errors of computed smallest and largest eigenpairs with dimensions $n=2,\dots,100$.}
\label{Ave_rel_berr}
\begin{tabular}{ c| c| c| c| c}
 &\multicolumn{2}{r|}{LINEARIZATION AVG.\ BK.\ ERR.}& \multicolumn{2}{c}{HOMOTOPY AVG.\ BK.\ ERR.} \\
\cline{2-5}
  $n$ &SMALLEST&LARGEST& SMALLEST&LARGEST  \\
\hline
2 & 2.02195E-16&1.07268E-16 & 1.53449E-16&1.38260E-16 \\
5 & 3.21357E-16&2.63883E-16 & 1.36922E-16&8.70572E-17	 \\
10 & 3.97923E-16&4.17780E-16 & 1.95933E-16&1.37749E-16	 \\
15 & 5.64010E-16&4.37375E-16 & 2.82767E-16&1.26614E-16	 \\
20 & 6.25007E-16&5.71766E-16 & 1.42798E-16&1.49327E-16	 \\
30 & 8.03800E-16&5.77483E-16 & 1.25286E-16&3.20068E-16	 \\
40 & 9.85766E-16&6.64878E-16 & 1.49583E-16&1.29559E-16	 \\
50 & 9.98415E-16&7.09502E-16 & 1.92299E-16&2.26032E-15	 \\
60 & 1.11238E-15&6.56447E-16 & 1.59075E-16&2.24618E-16	 \\
70 & 1.21839E-15&6.95665E-16 & 1.63888E-16&2.10081E-16	 \\
80 & 1.29842E-15&8.85652E-16 & 1.65893E-16&2.42794E-16	 \\
90 & 1.30478E-15&8.85652E-16 &  1.72905E-16&1.66882E-16	 \\
100 & 1.44644E-15&8.90062E-16 &	1.80208E-16&1.84720E-16
\end{tabular}
\end{table}

\subsection{Speed and accuracy comparisons for PEP}\label{PEP_test}

\begin{figure}[H]
\caption{Graphs for \cref{Ave_abs_berr_PEP}.
Blue for linearization and red for homotopy.}
\centering
\includegraphics[scale=0.36]{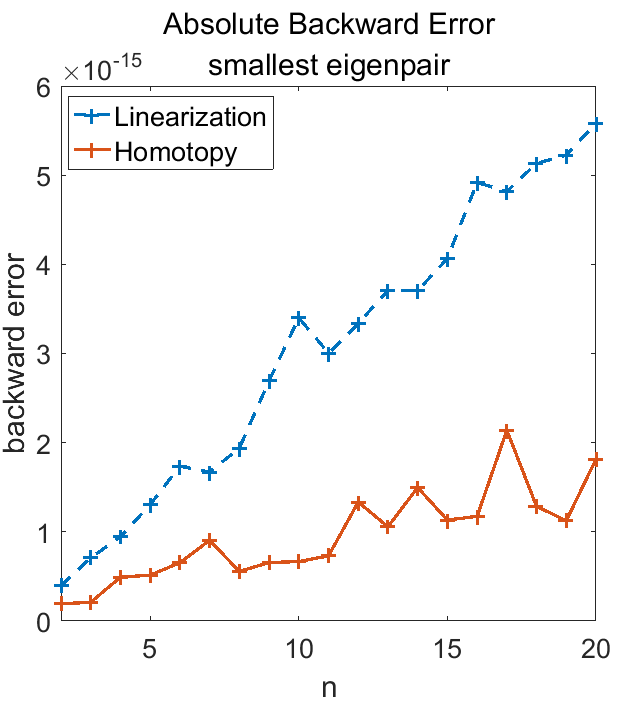}
\includegraphics[scale=0.36]{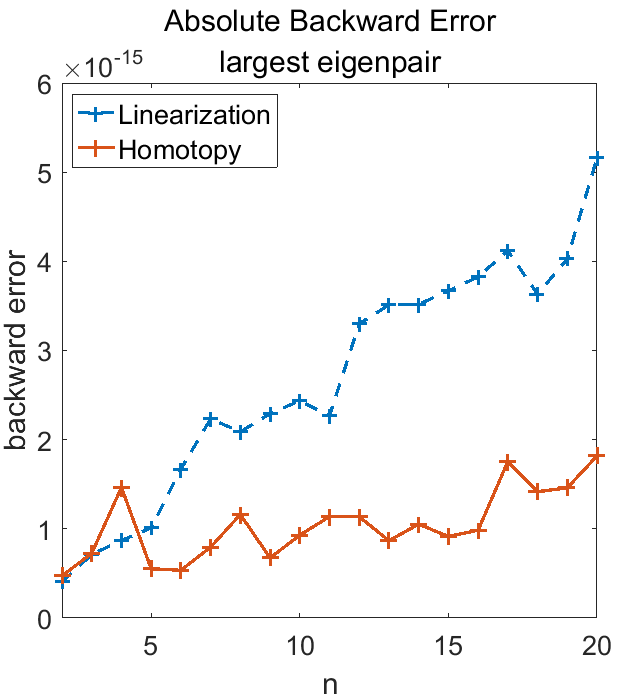}\\
\includegraphics[scale=0.36]{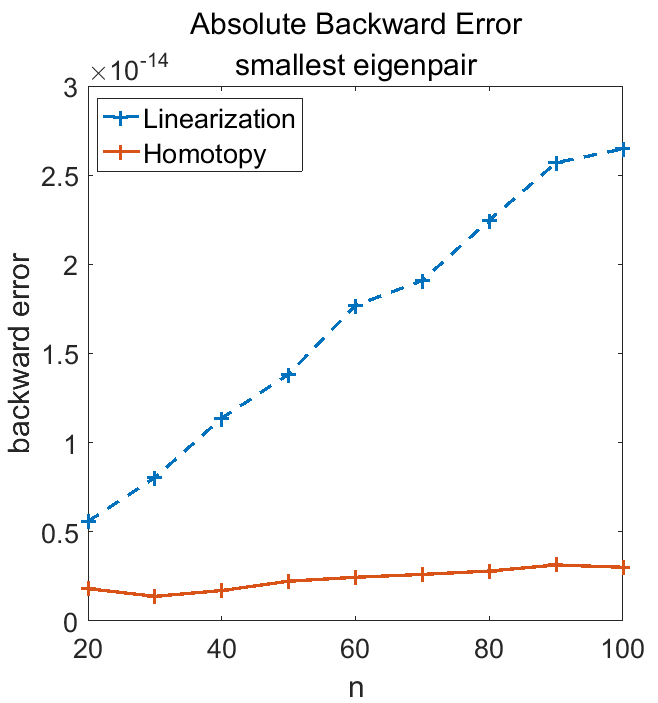}
\includegraphics[scale=0.36]{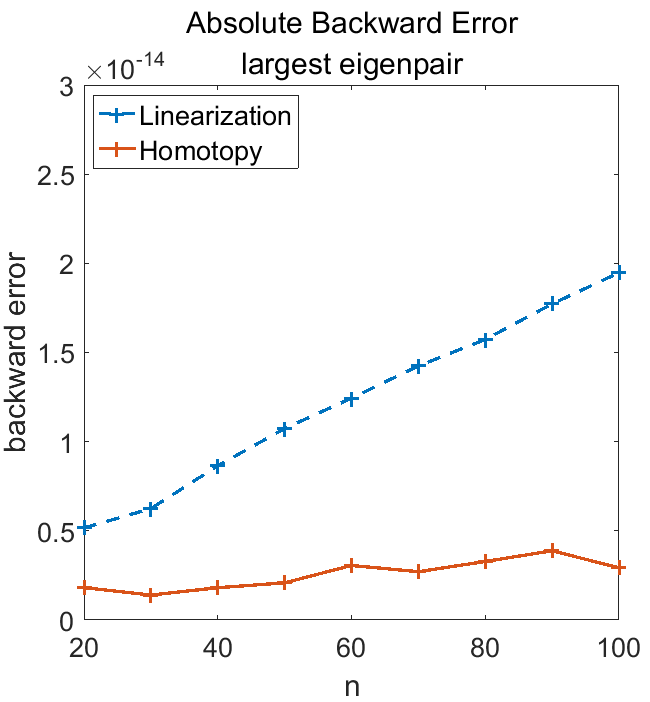}
\label{Abs_BErr_pep}
\end{figure}

\begin{figure}[H]
\caption{Graphs for \cref{Ave_rel_berr_PEP}.
Blue for linearization and red for homotopy.}
\centering
\includegraphics[scale=0.36]{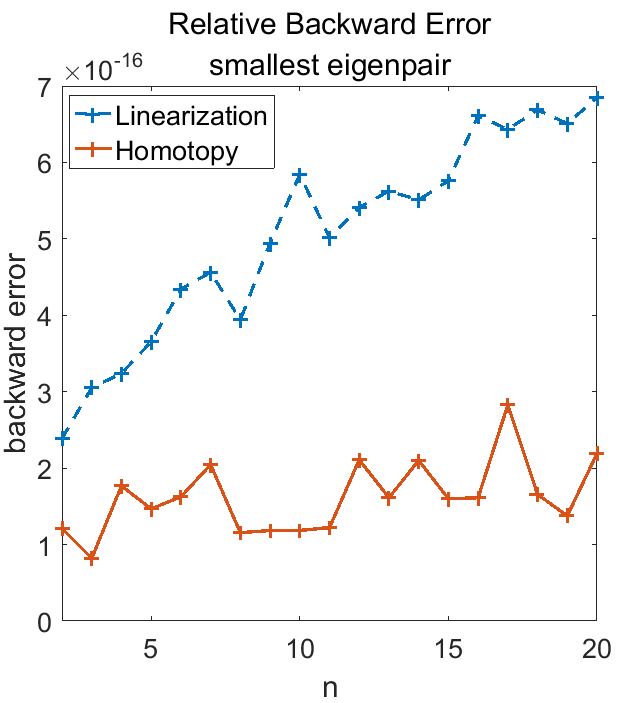}
\includegraphics[scale=0.36]{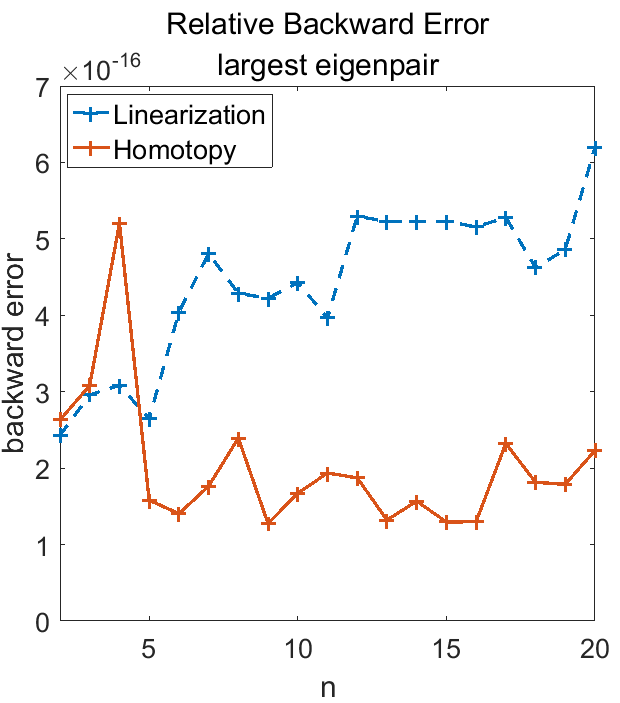}\\
\includegraphics[scale=0.36]{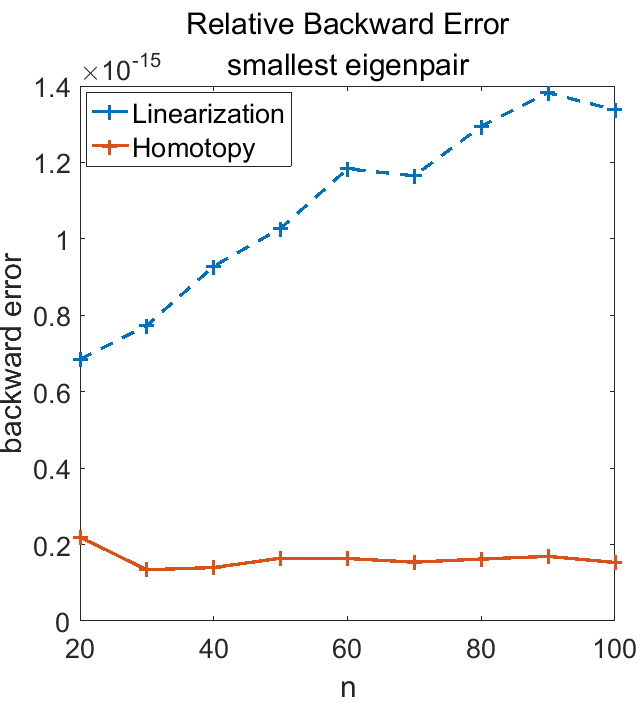}
\includegraphics[scale=0.36]{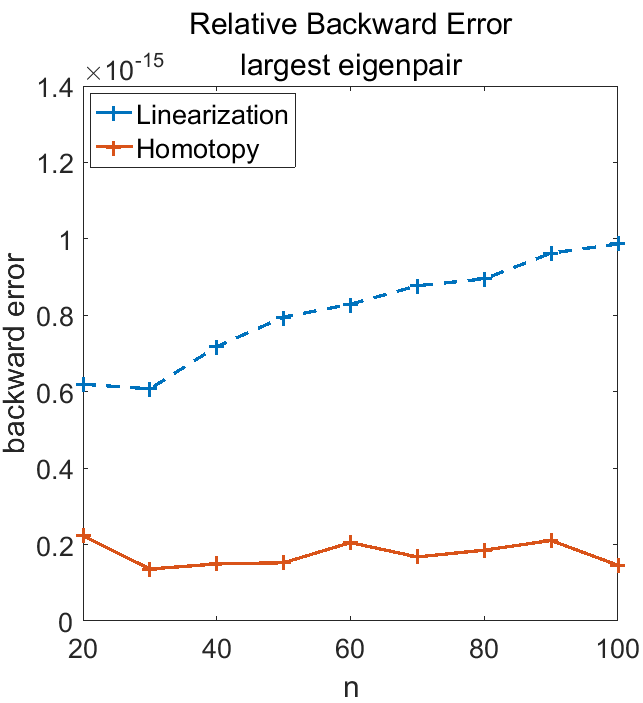}
\label{Rel_BErr_pep}
\end{figure}

We repeat the timing and accuracy tests in \cref{sec:speedQEP} and \cref{sec:accuracyQEP} for PEP where $m=4$. The timing results are presented in \cref{timing_linearization_PEP} and \cref{timing_homotopy_PEP}. The accuracy results are presented side-by-side in tabulated form in \cref{Ave_abs_berr_PEP} (absolute error) and \cref{Ave_rel_berr_PEP} (relative error), and graphically in \cref{Abs_BErr_pep} (absolute error) and \cref{Rel_BErr_pep} (relative error). As in the case of QEP,  for speed comparisons, we run our experiments ten times and record the best, average, median, and worst performance for each method and each dimension; for accuracy comparisons, the tests are averaged over ten runs. The conclusions drawn from these experiments are discussed in \cref{sec:conclude}.
\begin{table}[H]
\footnotesize
\centering
\caption{\textsc{Speed of PEP --- linearization.}
Elapsed timings (in seconds) for the linearization method with dimension $n=2,3,\dots,100$  (not all displayed).}
\label{timing_linearization_PEP}
\begin{tabular}{ c| c| c| c| c|c}
  $n$ & \# ROOTS &  BEST & AVERAGE & MEDIAN & WORST \\
\hline
2 & 8 & 3.6209E-4&	7.3132E-4	&6.2914E-4&	0.0016 \\
4 & 16 & 6.4617E-4&	9.8930E-4	&8.6807E-4&	0.0017 \\
6 & 24 & 0.0010&	0.0017	&0.0015&	0.0028 \\
8 & 32 & 0.0017&	0.0022	&0.0022&	0.0031 \\
10 & 40 & 0.0022&	0.0030	&0.0031&	0.0042 \\
12 & 48 & 0.0030&	0.0040	&0.0037&	0.0055 \\
14 & 56 & 0.0045&	0.0054	&0.0050&	0.0074 \\
16 & 64 & 0.0068&	0.0080	&0.0080&	0.0100 \\
18 & 72 & 0.0091&	0.0108	&0.0099&	0.0141 \\
20 & 80 & 0.0106&	0.0121	&0.0116&	0.0149 \\
30 & 120 & 0.0262&	0.0337	&0.0349&	0.0373 \\
40 & 160 & 0.0656&	0.0751	&0.0695&	0.0929 \\
50 & 200 & 0.0920&	0.1225	&0.1095&	0.1724 \\
60 & 240 & 0.1606&	0.1906	&0.1807&	0.2513 \\
70 & 280 & 0.2401&	0.2719	&0.2675&	0.3262 \\
80 & 320 & 0.3757&	0.4189	&0.4140&	0.4986 \\
90 & 360 & 0.4947&	0.5507	&0.5296&	0.6968 \\
100 & 400 & 0.7903&	0.8285	&0.8210&	0.9458
\end{tabular}
\end{table}

\begin{table}[H]
\footnotesize
\centering
\caption{\textsc{Speed of PEP --- homotopy in serial and in parallel.}
Elapsed timings (in seconds) for the homotopy method with dimension $n=2,3,\dots,100$; in serial for $n = 2,3,\dots,20$  (not all displayed); in parallel on $20$ cores for $n = 30,40,\dots,100$.}
\label{timing_homotopy_PEP}
\begin{tabular}{ c| c| c| c| c|c}
  n & \# ROOTS &  BEST & AVERAGE & MEDIAN & WORST \\
\hline
2 & 8 & 0.4449&	0.5051	&0.4725&	0.6503 \\
4 & 16 & 0.6237&0.8137	&0.6946&	1.2260 \\
6 & 24 & 1.0655&2.4848	&2.0590&4.9239 \\
8 & 32 & 1.6349&3.5440	&2.7851&7.6895 \\
10 & 40 &2.7746 &8.6236	&9.2868&16.6303 \\
12 & 48 &7.3469& 16.3581&16.8977&	25.2626 \\
14 & 56 & 10.0197&	26.8512	&18.3648&57.4097 \\
16 & 64 & 22.8207&	45.5760	&46.3289&66.0404 \\
18 & 72 & 19.6662&	61.4427	&56.9810&94.0309 \\
20 & 80 & 47.8803&	119.9619	&98.6018&	245.0955 \\
\hline
30 & 120 & 14.0910&	30.9721	&30.2710& 56.8700 \\
40 & 160 & 35.1230&	55.5960	&51.7820& 77.9600 \\
50 & 200 & 75.2800&	111.3593	&111.8870&	148.0980 \\
60 & 240 & 165.5460&214.3197	&205.0230&	290.4680 \\
70 & 280 & 290.4960&397.6789	&421.8780&	474.0510 \\
80 & 320 & 479.8200&616.7326	&574.7630&	811.5990 \\
90 & 360 & 684.9450&1037.8848	&1114.4720&	1258.1440 \\
100 & 400 & 1046.2700&	1452.7988	&1423.2730&	1947.1600
\end{tabular}
\end{table}

\begin{table}[H]
\footnotesize
\centering
\caption{\textsc{Accuracy of PEP --- linearization vs homotopy.}
Absolute backward errors of computed smallest and largest eigenpairs with dimension $n=2,\dots,100$.}
\label{Ave_abs_berr_PEP}
\begin{tabular}{ c| c| c| c| c}
 &\multicolumn{2}{r|}{LINEARIZATION AVG.\ BK.\ ERR.}& \multicolumn{2}{c}{HOMOTOPY AVG.\ BK.\ ERR.} \\
\cline{2-5}
  $n$ &SMALLEST&LARGEST& SMALLEST&LARGEST  \\
\hline
2 & 3.97673E-16&4.06234E-16 &1.93084E-16 & 4.72180E-16\\
5 & 1.30486E-15&1.01101E-15 &5.14466E-16 & 5.52493E-16\\
10 & 3.40010E-15&2.43585E-15 &6.66926E-16 &9.29498E-16\\
15 & 4.06465E-15&3.66940E-15 &1.13208E-15 &9.14666E-16\\
20 & 5.58414E-15&5.16130E-15 &1.81496E-15 &1.81759E-15\\
30 & 8.03892E-15&6.24574E-15 &1.38205E-15 &1.39749E-15\\
40 & 1.13754E-14&8.64433E-15 &1.70770E-15 &1.81052E-15\\
50 & 1.38394E-14&1.07375E-14 &2.23095E-15 &2.08415E-15\\
60 & 1.76932E-14&1.24398E-14 &2.45192E-15 &3.05775E-15\\
70 & 1.90902E-14&1.42487E-14 &2.61274E-15 &2.70850E-15\\
80 & 2.24766E-14&1.57303E-14 &2.80187E-15 &3.28150E-15\\
90 & 2.56978E-14&1.77254E-14 &3.14128E-15 &3.88126E-15\\
100 & 2.64759E-14&1.94732E-14&3.01357E-15 &2.92889E-15
\end{tabular}
\end{table}

\begin{table}[H]
\footnotesize
\centering
\caption{\textsc{Accuracy of PEP --- linearization vs homotopy.}
Relative backward errors of computed smallest and largest eigenpairs with dimension $n=2,\dots,100$.}
\label{Ave_rel_berr_PEP}
\begin{tabular}{ c| c| c| c| c}
 &\multicolumn{2}{r|}{LINEARIZATION AVG.\ BK.\ ERR.}& \multicolumn{2}{c}{HOMOTOPY AVG.\ BK.\ ERR.} \\
\cline{2-5}
  $n$ &SMALLEST&LARGEST& SMALLEST&LARGEST  \\
\hline
2 & 2.38800E-16&2.42932E-16 &1.20966E-16 &2.63005E-16\\
5 & 3.65363E-16&2.64597E-16 &1.46731E-16 &1.58114E-16\\
10 & 5.84645E-16&4.42737E-16 &1.18572E-16 &1.67121E-16\\
15 & 5.76099E-16&5.22943E-16 &1.60022E-16 &1.29557E-16\\
20 & 6.84554E-16&6.19752E-16 &2.19221E-16 &2.22975E-16\\
30 & 7.73547E-16&6.07382E-16 &1.33975E-16 &1.36049E-16\\
40 & 9.28358E-16&7.18160E-16 &1.39950E-16 &1.49331E-16\\
50 & 1.02753E-15&7.95390E-16  &1.64387E-16 &1.52488E-16\\
60 & 1.18312E-15&8.28990E-16 &1.63573E-16 &2.05105E-16\\
70 & 1.16505E-15&8.77305E-16 &1.53855E-16 &1.67640E-16\\
80 & 1.29478E-15&8.95204E-16 &1.62206E-16 &1.85393E-16\\
90 & 1.38225E-15&9.62668E-16 &1.69063E-16 &2.10661E-16\\
100 & 1.33725E-15&9.86613E-16&1.53193E-16 &1.45843E-16
\end{tabular}
\end{table}

\subsection{Conclusions}\label{sec:conclude}

From the results of the numerical experiments in preceding subsections, we may draw several conclusions regarding the speed and accuracy of our proposed homotopy method versus linearization method, the existing method-of-choice.
\begin{description}
\item[Speed] From the elapsed times, the linearization method is consistently much faster than the homotopy method, even when the latter is run in parallel.
\item[Accuracy] From the  normwise backward errors, homotopy method is more accurate than the linearization method across coefficient matrices of all dimensions.
\item[Dimension] The gap in accuracy increases significantly  with the dimension of the coefficient matrices, favoring homotopy method for higher dimensional problems, particularly when all eigenpairs are needed.
\item[Stability] Our results reflect what is known about linearization method (see \cref{ss:Accuracy}), namely, it is numerically unstable for both QEP and PEP.
\item[Initialization] Our results reflect what is known about homotopy method (see \cref{ss:homotopymethod}), namely, its dependence on the choice of start system --- one that leads to a path passing near a singularity takes longer to converge.
\item[Parallelism] The timings of homotopy method can be improved substantially with parallel computing. On the other hand, the linearization method exhibits no obvious parallelism.
\end{description}

It is also evident from these results that both the linearization and homotopy methods take considerably longer to solve a PEP with $m=4$ than to solve a QEP. This can be attributed to the fact that a PEP with $m =4$ has $5/3$ times more parameters than a QEP: $\operatorname{dim} (\mathbb{C}^{n \times n \times 5}) = 5n^2$ versus $\operatorname{dim} (\mathbb{C}^{n \times n \times 3}) = 3n^2$, and that a PEP with $m=4$ has twice as many eigenpairs as a QEP: $4n$ versus $2n$ eigenpairs.

It is interesting to observe that accuracy does not exhibit a similar deterioration with increasing $m$ --- for a fixed $n$, the averaged backward error for a PEP with $m=4$ can be smaller than that for a QEP. 

\section{Conditioning and accuracy}\label{sec:cond} 

It has been shown \cite{doi:10.1137/050628283} that if the $2$-norms of the coefficient matrices in a PEP are all approximately $1$, then the companion linearization and original PEP have similar condition numbers. In particular, define
\begin{equation}\label{rho}
\rho = \frac{\max_i \|A_i \|_2}{\min(\|A_0 \|_2,\|A_m\|_2)} \geq 1.
\end{equation}
When $\rho$ is of order 1, there exists a linearization for a particular eigenvalue that is about as well conditioned as the original PEP itself for that eigenvalue, to within a small constant factor. However this is no longer true when the $A_i$'s vary widely in norm: The companion GEP is potentially far more ill-conditioned than the PEP. 

For the QEP in \cref{QEP_def}, the quantity
\[
\rho = \frac{\max(\|M\|,\|C\|,\|K\|)}{\min(\|M\|,\|K\|)}
\]
is of order $1$ if
$\|C\| \lesssim \max(\|M\|,\|K\|)$ and $\|M\| \approx \|K\|$.
When these are not satisfied, a scaling of $Q(\lambda)$ will typically improve the conditioning of the linearization --- provided that $Q(\lambda)$ is not too heavily damped, i.e., $\|C\|\lesssim \sqrt{\|M\|\|K\|}$.
However, in general such a scaling is unavailable; for instance, it is still not known how one should scale a heavily damped QEP.

\begin{figure}[H]
\caption{Relative backward errors (in Frobenius norm) of all computed eigenpairs.  Blue dots: linearization method. Red crosses: homotopy method.}
\centering
\includegraphics[scale=0.41]{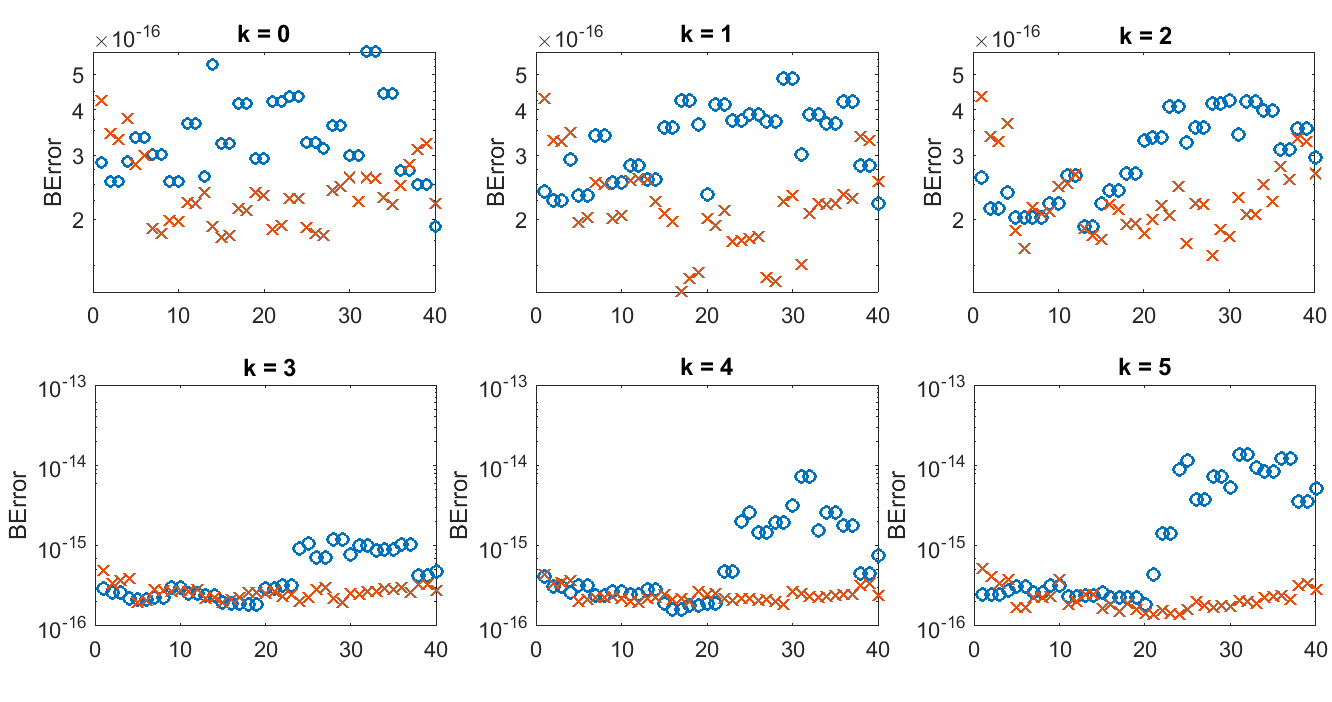}\\
\label{RBError_ScaleC}
\end{figure}

\begin{figure}[H]
\caption{Condition numbers of all computed eigenpairs. Blue dots: companion GEPs. Red crosses: original QEPs.}
\centering
\includegraphics[scale=0.41]{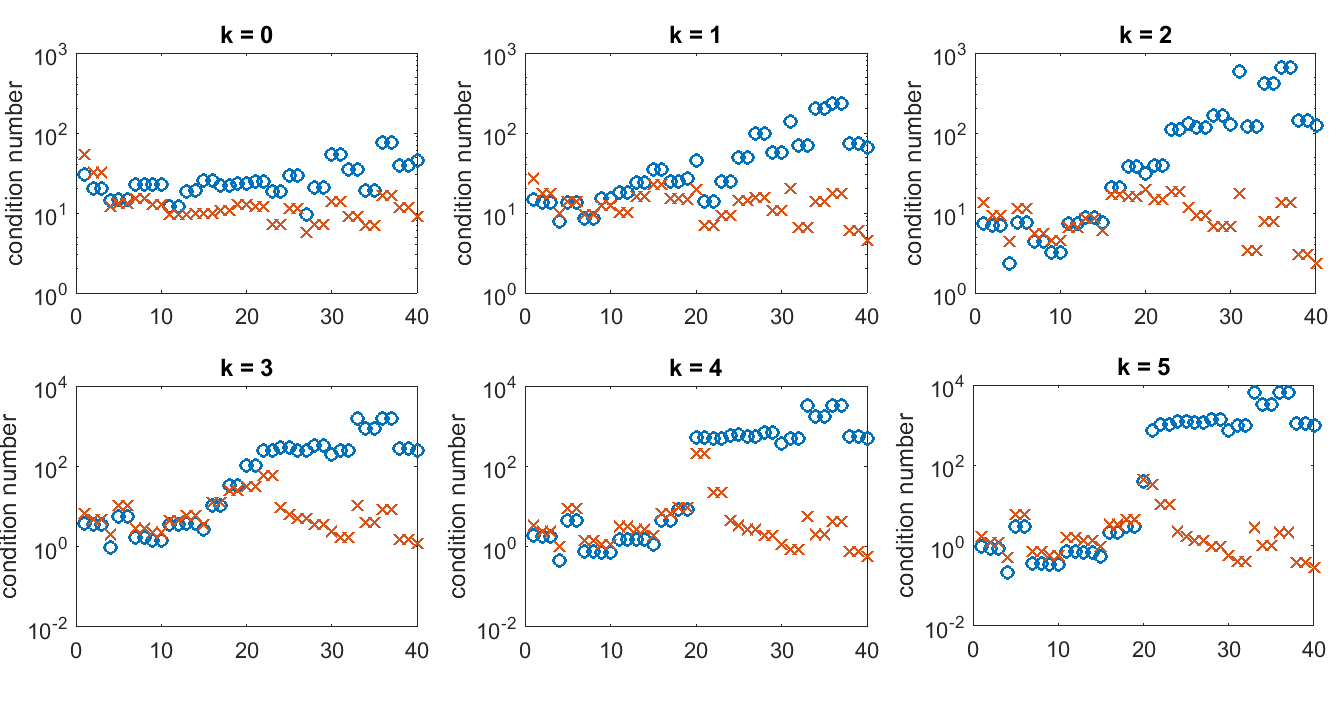}
\label{ConditionNumber_ScaleC}
\end{figure}

We compare how linearization and homotopy methods perform on damped QEPs:
\begin{enumerate}[\upshape (i)]
\item We generate random $20 \times 20$ coefficient matrices $M,C,K$ having independent entries that follow the standard real Gaussian distribution $\mathcal{N}_\mathbb{R}(0,1)$.
\item For $Q_k(\lambda)=\lambda^2M+\lambda (2^k\cdot C)+K$, $k=0,1,\dots,5$, we determine the relative backward errors of all computed eigenpairs for both linearization and homotopy methods. The results are in \cref{RBError_ScaleC}.
\item At the same time, for each $Q_k(\lambda)$, we compute the condition number of each eigenvalue  in both the original QEP (used in homotopy method) and its companion GEP (used in linearization method). The results are in \cref{ConditionNumber_ScaleC}.
\end{enumerate}
In both \cref{RBError_ScaleC} and  \cref{ConditionNumber_ScaleC}, the horizontal axis is the index of eigenvalues in ascending order of magnitude, the vertical axis is on a log scale.

We may deduce the following from these results:
\begin{enumerate}[\upshape (a)]
\item From \cref{RBError_ScaleC} we see that homotopy method is backward stable for all eigenvalues and all $k = 0, 1, \dots, 5$; on the other hand, linearization method becomes significantly less stable for the larger eigenvalues as $k$ increases past $3$. Note that the larger the value of $k$, the heavier damped the QEP.
\item From \cref{ConditionNumber_ScaleC} we see that the larger eigenvalues in the original QEP are far better conditioned than in its companion GEP, whereas the smaller eigenvalues are similarly conditioned in both problems.
\item From \cref{ConditionNumber_ScaleC} we see that as the QEP becomes more heavily damped, the companion GEP becomes exceedingly worse-conditioned than the original QEP.
\end{enumerate}
In summary, homotopy method is evidently the preferred method for accurate determination of large eigenvalues in heavily damped QEPs.

Although it is possible to find an alternative linearization of a QEP into a GEP that is better conditioned towards the larger eigenvalues \cite{tisseur2000backward}, there is no known linearization with conditioning comparable to the original QEP across \emph{all} eigenvalues (in fact, such an ideal linearization quite likely does not exist). In principle, one might use different linearizations to determine eigenvalues in different ranges, and then combine the results to obtain a full set of eigenvalues. However this is not only impractical but suffers from a fallacy --- we do not know a priori which eigenvalues from which linearizations are more accurate.

\section{Applications}\label{sec:app}

In this section we provide numerical experiments on real data (as opposed to randomly generated data in \cref{Numerical results}) arising from the two application problems mentioned in \cref{Intro}. All experiments here are conducted in the same environment as described in \cref{Numerical results}, except that the matrices defining the problems are generated with the \texttt{nlevp} function in the \textsc{Matlab} toolbox \textsc{Nlevp} \cite{MIMS_ep2011.116}.

\subsection{Acoustic wave problem} \label{ss:QEPEx}

This application is taken from \cite{acoustic_wave}. Consider an acoustic medium with constant density and space-varying sound speed $c(x)$ occupying the volume $\Omega \subseteq  \mathbb{R}^d$. 
The homogeneous wave equation for the acoustic pressure $p(x,t)=\widehat p(x)e^{\widehat{\lambda} t}$
has a factored form that simplifies the wave equation to the following, where  $\widehat{\lambda}$ is the eigenvalue:
\begin{equation}\label{time_harmonic}
-\Delta \widehat{p}(x)-\bigl(\widehat{\lambda}/c\bigr)^2 \cdot \widehat{p}(x)=0.
\end{equation}
For our purpose, it suffices to consider the one-dimensional acoustic wave problem, i.e., $d=1$ and $\Omega=[0,1]$, with Dirichlet  boundary condition $\widehat{p} = 0$ and impedance boundary condition $\partial p / \partial n + i\widehat\lambda p/ \zeta = 0$.

\subsubsection{Numerical results} The quadratic matrix polynomial $Q(\lambda)=\lambda^2 M+ \lambda C+K$ arises from a finite element discretization of \cref{time_harmonic}.
In our numerical experiments we set  impedance $\zeta=1$. The three $n\times n$ matrices in this QEP are
\[
M=-4\pi^2 \frac{1}{n}\left(I_n-\frac{1}{2}e_ne_n^\tp\right),
\quad C=2\pi i \frac{1}{\zeta}e_ne_n^\tp,
\quad K=
n\begin{bmatrix}      2 & -1 & & \\     -1 &\ddots & \ddots &\\	& \ddots &2&-1\\    & & -1&1     \end{bmatrix}.
\]

\begin{figure}[H]
\caption{\textsc{Acoustic wave problem --- absolute and relative backward errors.}  Computed smallest and largest eigenpairs with dimensions $n=20,30,\dots,100$.
Blue dashed lines: linearization method. Red solid lines: homotopy method.}
\centering
\includegraphics[scale=0.36]{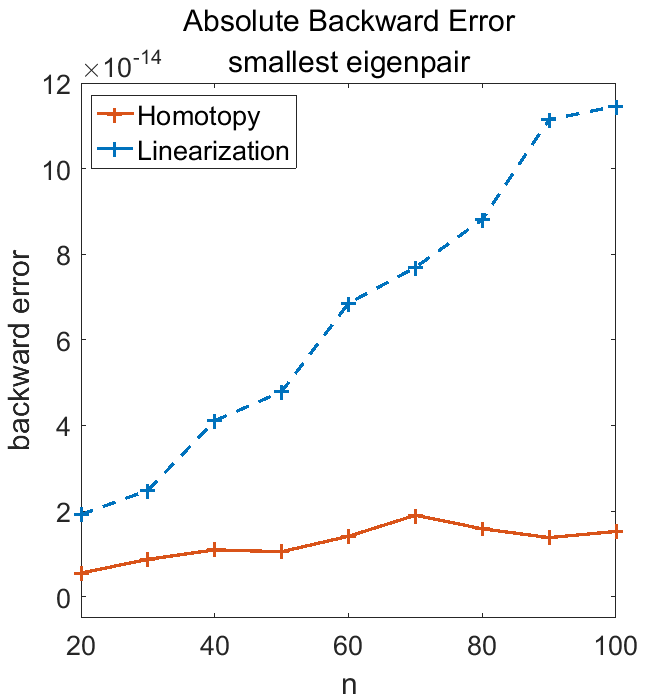}
\includegraphics[scale=0.36]{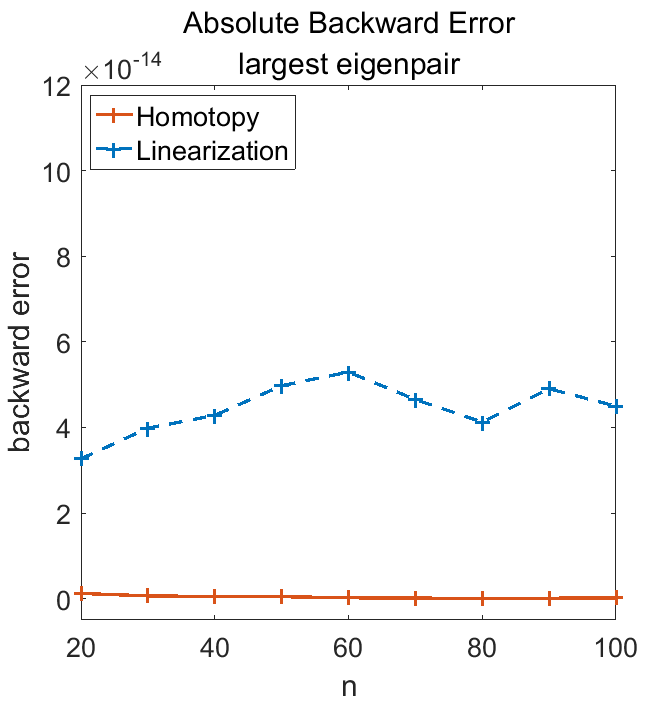}\\
\includegraphics[scale=0.36]{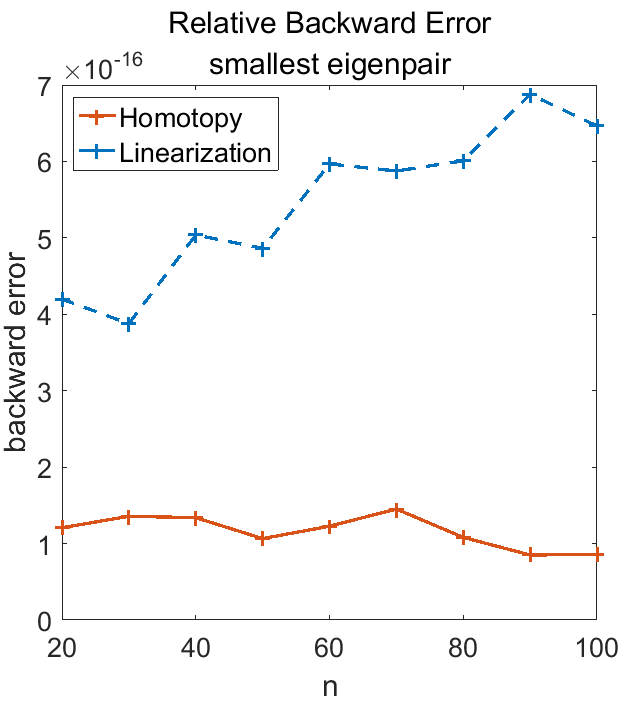}
\includegraphics[scale=0.36]{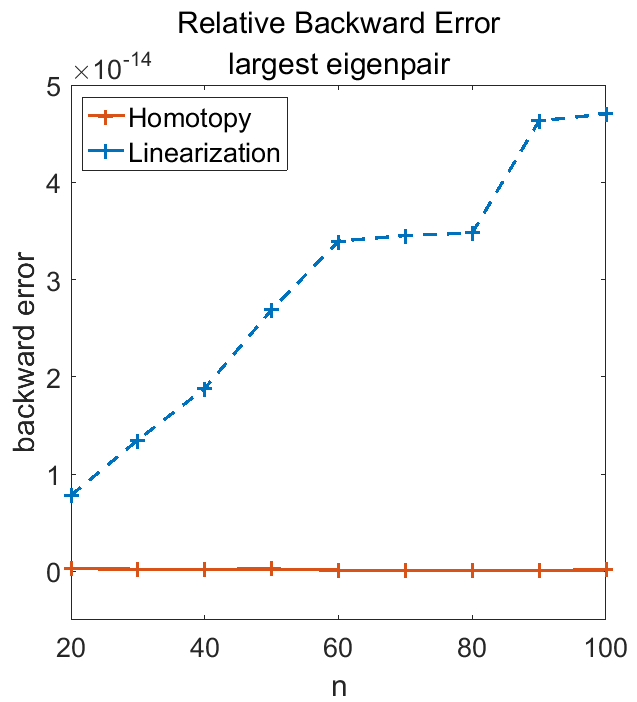}
\label{ex1_BErr_20-100}
\end{figure}

As before, we compare the accuracy and timings of homotopy and linearization methods on this problem. We tabulate the absolute and relative backward errors of the computed eigenpairs corresponding to the smallest and largest eigenvalues  for  dimensions $n=20,30,\dots,100$  in \cref{Ave_abs_berr_APP1} and \cref{Ave_rel_berr_APP1}, and plot them graphically in \cref{ex1_BErr_20-100}.  All accuracy tests are averaged over ten runs.
The elapsed timings are tabulated in \cref{timing_linearization_APP1} and \cref{timing_homotopy_APP1}.  All speed tests are run ten times with  the best, average, median and worst timings recorded. 

\begin{table}[H]
\footnotesize
\caption{\textsc{Acoustic wave problem ---
 absolute backward errors.} Computed smallest and largest eigenpairs with dimensions $n=20,30,\dots,100$ via homotopy and linearization methods.}
\label{Ave_abs_berr_APP1}
\begin{center}
\begin{tabular}{ c| c| c| c| c}
 &\multicolumn{2}{r|}{LINEARIZATION AVG.\ BK.\ ERR.}& \multicolumn{2}{c}{HOMOTOPY AVG.\ BK.\ ERR.} \\
\cline{2-5}
  n &SMALLEST&LARGEST& SMALLEST&LARGEST  \\
\hline
20 & 1.91540E-14&3.26558E-14 & 5.50499E-15&1.20102E-15 \\
30 & 2.48304E-14&3.98265E-14 & 8.69382E-15&6.41943E-16	 \\
40 & 4.10844E-14&4.28048E-14 & 1.09178E-14&4.76802E-16	 \\
50 & 4.78653E-14&4.97587E-14 & 1.04868E-14&4.50403E-16	 \\
60 & 6.85333E-14&5.29234E-14 & 1.40925E-14&1.91169E-16	 \\
70 & 7.68982E-14&4.65017E-14 & 1.90018E-14&1.07398E-16	 \\
80 & 8.80830E-14&4.12431E-14 & 1.58457E-14&8.02790E-17	 \\
90 & 1.11451E-13&4.90212E-14 & 1.37846E-14&9.11717E-17	 \\
100 & 1.14562E-13&4.49546E-14 & 1.52145E-14&1.52150E-16	
\end{tabular}
\end{center}
\end{table}

\begin{table}[H]
\footnotesize
\caption{\textsc{Acoustic wave problem ---
relative backward errors.}  Computed smallest and largest eigenpairs with dimensions $n=20,30,\dots,100$ via homotopy and linearization methods.}
\label{Ave_rel_berr_APP1}
\begin{center}
\begin{tabular}{ c| c| c| c| c}
 &\multicolumn{2}{r|}{LINEARIZATION AVG.\ BK.\ ERR.}& \multicolumn{2}{c}{HOMOTOPY AVG.\ BK.\ ERR.} \\
\cline{2-5}
  n &SMALLEST&LARGEST& SMALLEST&LARGEST  \\
\hline
20 & 4.19838E-16&7.81869E-15 & 1.20664E-16&2.87557E-16 \\
30 & 3.87748E-16&1.35042E-14 & 1.35761E-16&2.17667E-16	 \\
40 & 5.03756E-16&1.88157E-14 & 1.33868E-16&2.09588E-16	 \\
50 & 4.86198E-16&2.68895E-14 & 1.06521E-16&2.43396E-16	 \\
60 & 5.96657E-16&3.39445E-14 & 1.22691E-16&1.22614E-16	 \\
70 & 5.87499E-16&3.45264E-14 & 1.45173E-16&7.97402E-17	 \\
80 & 6.00842E-16&3.47937E-14 & 1.08089E-16&6.77252E-17	 \\
90 & 6.87825E-16&4.63156E-14 & 8.50724E-17&8.61397E-17	 \\
100 & 6.46406E-16&4.70234E-14 & 8.58465E-17&1.59152E-16	
\end{tabular}
\end{center}
\end{table}
\begin{table}[H]
\footnotesize
\centering
\caption{\textsc{Acoustic wave problem --- linearization method.}
Elapsed timings (in seconds) for the linearization method with dimensions $n=20, 30, \dots,100$.}
\label{timing_linearization_APP1}
\begin{tabular}{ c| c| c| c| c|c}
  n & \# ROOTS &  BEST & AVERAGE & MEDIAN & WORST \\
\hline
20 & 40 & 0.0056&	0.0139	&0.0071&	0.0754 \\
30 & 60 & 0.0107&	0.0175	&0.0134&	0.0592 \\
40 & 80 & 0.0202&	0.0310	&0.0221&	0.1127 \\
50 & 100 & 0.0365&	0.0499	&0.0387&	0.1464 \\
60 & 120 & 0.0469&	0.0673	&0.0600&	0.1380 \\
70 & 140 & 0.0680&	0.0826	&0.0830&	0.1040 \\
80 & 160 & 0.0933&	0.1192	&0.1004&	0.2491 \\
90 & 180 & 0.1304&	0.1549	&0.1445&	0.2135 \\
100 & 200 & 0.1723&	0.2015	&0.1818&	0.3137
\end{tabular}
\end{table}

\begin{table}[H]
\footnotesize
\centering
\caption{\textsc{Acoustic wave problem --- homotopy method.}
Elapsed timings (in seconds) for the homotopy method with dimensions $n=20, 30, \dots,100$, in serial for $n = 20,30,40$, and in parallel on $20$ cores for $n = 50,60,\dots,100$.}
\label{timing_homotopy_APP1}
\begin{tabular}{ c| c| c| c| c|c}
  n & \# ROOTS &  BEST & AVERAGE & MEDIAN & WORST \\
\hline
20 & 40 & 7.7472&	19.5057	&16.7535&	34.7300 \\
30 & 60 & 104.2553&	166.2379	&155.9914&	257.1145 \\
40 & 80 & 265.2321&	557.2342	&526.8778&	1003.7452 \\
\hline
50 & 100 & 30.0660&	37.9267	&35.5240&	49.5630 \\
60 & 120 & 74.2800&	92.1235	&94.7845&	113.7100 \\
70 & 140 & 149.3000&	189.5564	&182.7230&	280.6460 \\
80 & 160 & 252.4920&	378.7414	&370.6855&	519.2240 \\
90 & 180 & 530.7180&	691.8141	&736.6715&	841.8600 \\
100 & 200 & 559.9340&	1079.0177	&1154.3995&	1464.4170
\end{tabular}
\end{table}

\subsubsection{Conclusions} From these results, we may draw essentially the same conclusions in \cref{sec:conclude}  that we made about our numerical experiments on randomly generated data. In particularly,  homotopy method vastly outperforms linearization method in accuracy, as measured by normwise backward errors. The difference in this case is that we have very sparse, highly structured matrices $M,C,K$.


\subsubsection{Conditioning and accuracy} 
\begin{figure}[H]
\caption{\textsc{Acoustic wave problem with $n =100$ --- accuracy and conditioning.}  Horizontal axes: index of eigenvalues in ascending order of magnitude. Vertical axis: log scale.  \textsc{Left plot:} relative backward errors of computed eigenpairs; blue dots for the linearization method, red crosses for the homotopy method. \textsc{Right plot:} condition numbers of each eigenvalue; blue dots for the companion GEP, red crosses for the original QEP.}
\centering
\includegraphics[scale=0.36]{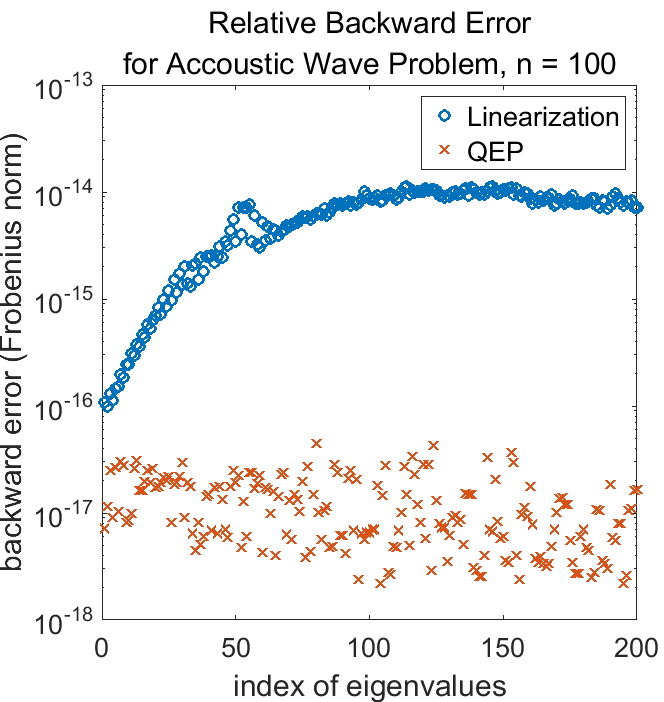}
\includegraphics[scale=0.37]{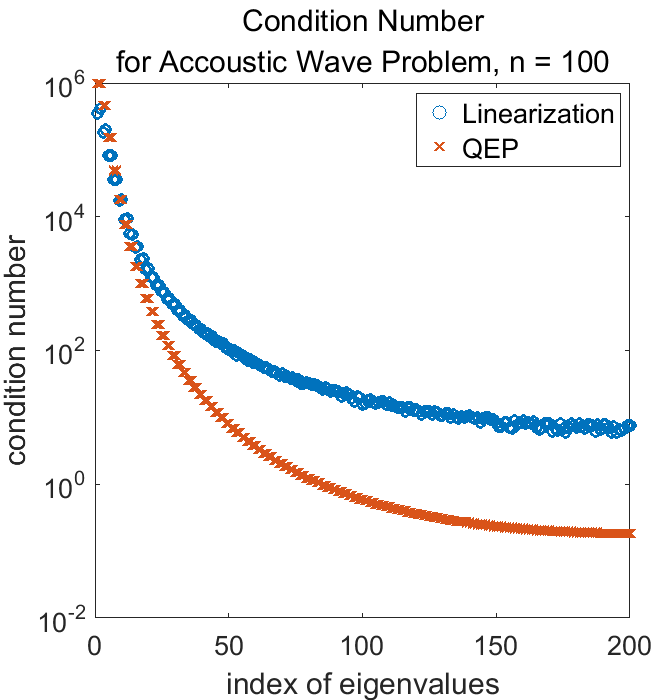}
\label{RBError_ConditionNumber_APP1}
\end{figure}

We next examine the accuracy of the two methods for the case $n=100$ more closely --- computing all eigenvalues, instead of just the largest and smallest. We will also compute the condition number for each eigenvalue in both the original QEP formulation and in the companion GEP formulation. The results are plotted in \cref{RBError_ConditionNumber_APP1}. In this case, $\|M\|_F=3.9$, $\|C\|_F=6.3$, $\|K\|_F=2439.3$, which is not regarded as a heavily damped QEP. 

Again, we see that homotopy method is numerically stable across all eigenpairs whereas the  linearization method is unstable for nearly 90\% of  eigenpairs. The original QEP is also better-conditioned than its companion GEP:  the condition number of each eigenvalue of the former is consistently smaller than, or very close to, that of the same eigenvalue in the latter.

\subsubsection{Certification} Last but not least, we will use Smale's $\alpha$-theory, discussed in \cref{sec:cert}, to deduce that for $n=20,60,100,140$, the Newton iterations in an application of homotopy method to the acoustic wave problem will converge quadratically to all eigenpairs.

\begin{figure}[H]
\caption{\textsc{Acoustic wave problem --- certification.} The values of $\beta(f,x,\lambda)$, $\mu(f,x,\lambda)$, and upper bound values of $\alpha(f,x,\lambda)$, $\gamma(f,x,\lambda)$ for the acoustic wave problem with $n=20,60,100,140$. 
The circle in the middle, the whiskers above and below represent the mean, maximum, and minimum value of each of these quantities taken over all eigenpairs. The red horizontal dashed line in the plot for $\alpha(f,x,\lambda)$ represents the threshold in \cref{smale}. Note that the vertical axis is in log scale.}
\centering
\includegraphics[scale=0.36]{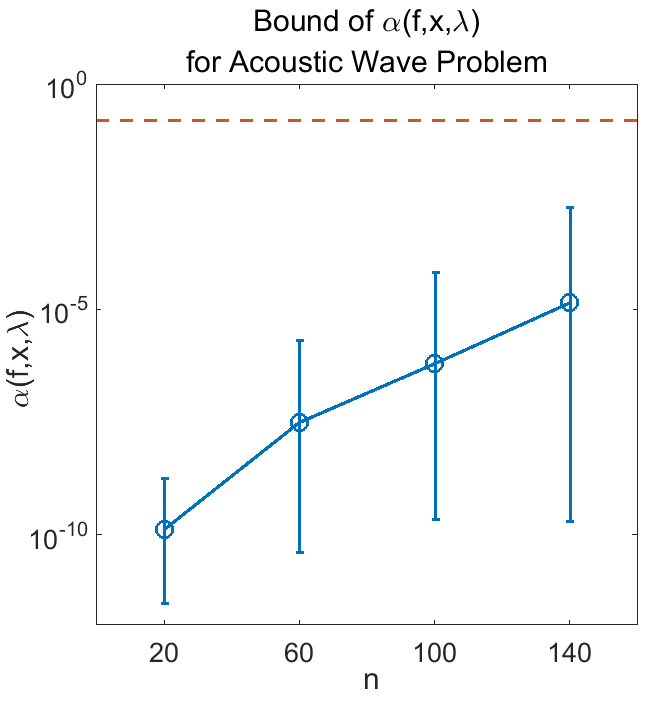}
\includegraphics[scale=0.36]{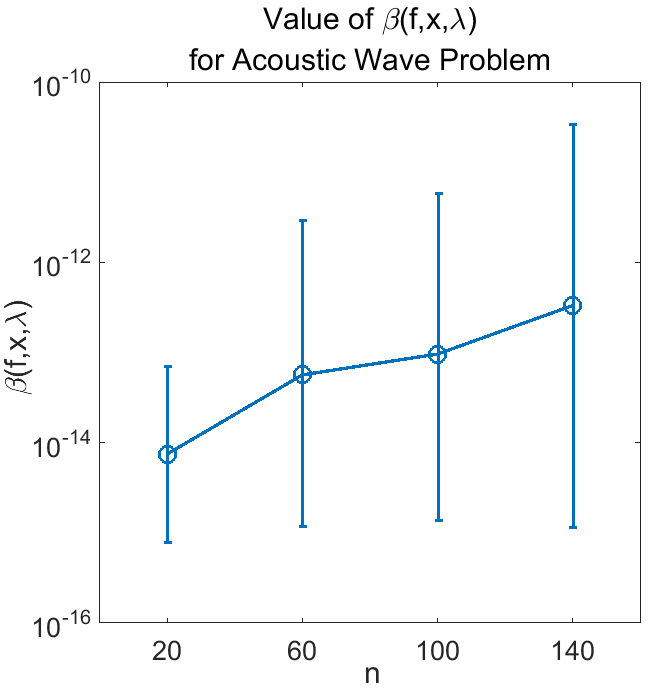}\\
\includegraphics[scale=0.36]{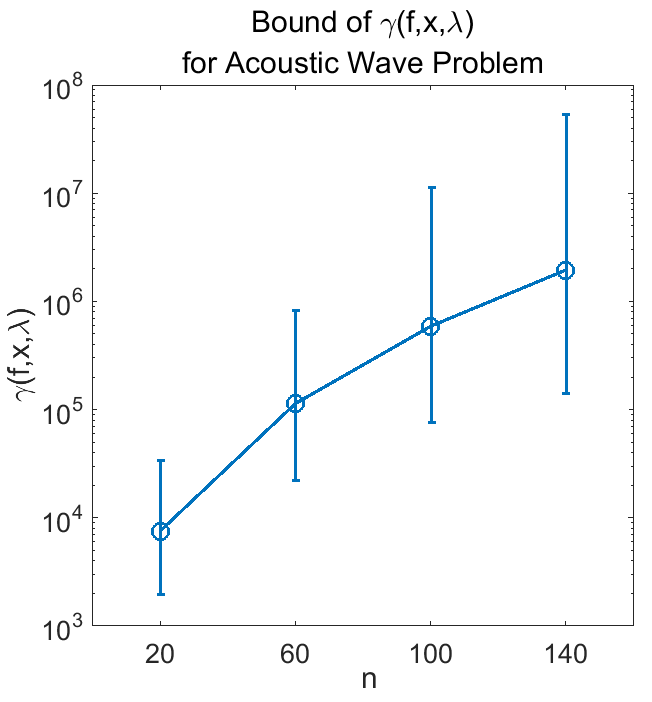}
\includegraphics[scale=0.36]{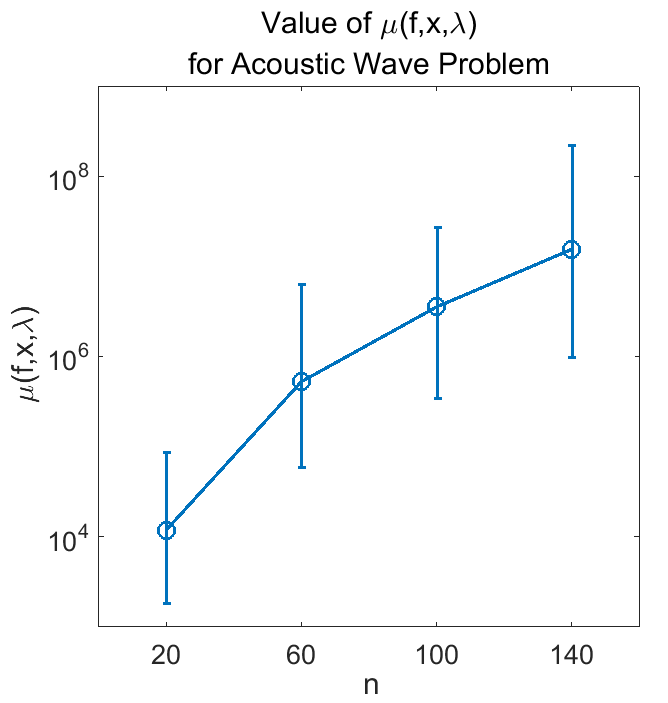}
\label{Certification_APP1}
\end{figure}

Specifically, for $n=20,60,100,140$, we solve for all eigenpairs, compute $\mu(f,x,\lambda)$ using \cref{mu_formula}, and obtain a bound for $\gamma(f,x,\lambda)$ according to \cref{pep_certification}. 
The Newton residual $\beta(f,x,\lambda)$ is also computed, which together with $\gamma(f,x,\lambda)$ yields  an upper bound for $\alpha(f,x,\lambda)$.  
We present our results in  \cref{Certification_APP1}.

For every dimension that we test,  the value of $\alpha(f,x,\lambda)$ for each eigenpair is much smaller than the threshold given in \cref{smale}. In other words, this certifies that all eigenpairs that we computed using the homotopy method are  accurate solutions  to the PEP,  for $n=20,60,100,140$.

\subsection{Planar waveguide problem} \label{ss:PEPEx}

This example is taken from \cite{stowell2010guided}. The $129\times 129$ quartic matrix polynomial $P(\lambda)=\lambda^4A_4+\lambda^3A_3+\lambda^2A_2+\lambda A_1+A_0$ arises from a finite element solution of the equation for the modes of a planar waveguide using piecewise linear basis $\varphi_i$, $i=0,\dots,128$. The coefficient matrices are defined by:
\begin{gather*}
A_1=\frac{\delta^2}{4}\operatorname{diag}(-1,0,0,\dots,0,0,1), \quad A_3=\operatorname{diag}(1,0,0,\dots,0,0,1),\\
A_0(i,j)=\frac{\delta^4}{16}\langle \varphi_i,\varphi_j\rangle,\quad A_2(i,j)=\langle\varphi_i',\varphi_j'\rangle -\langle q\varphi_i,\varphi_j\rangle,\quad A_4(i,j)=\langle\varphi_i,\varphi_j\rangle.
\end{gather*}
The parameter $\delta$ describes the difference in refractive index between the cover and  substrate of the waveguide and $q$ is a function from the variational formulation.

\begin{figure}[H]
\caption{\textsc{Planar waveguide problem --- sorted absolute and relative backward errors.}
Sorted absolute and relative backward errors of all $516$ computed eigenpairs. Blue dashed lines represent linearization method; red solid lines represent homotopy method. Tests are averaged over ten runs and vertical axis is in log scale.}
\centering
\includegraphics[scale=0.43]{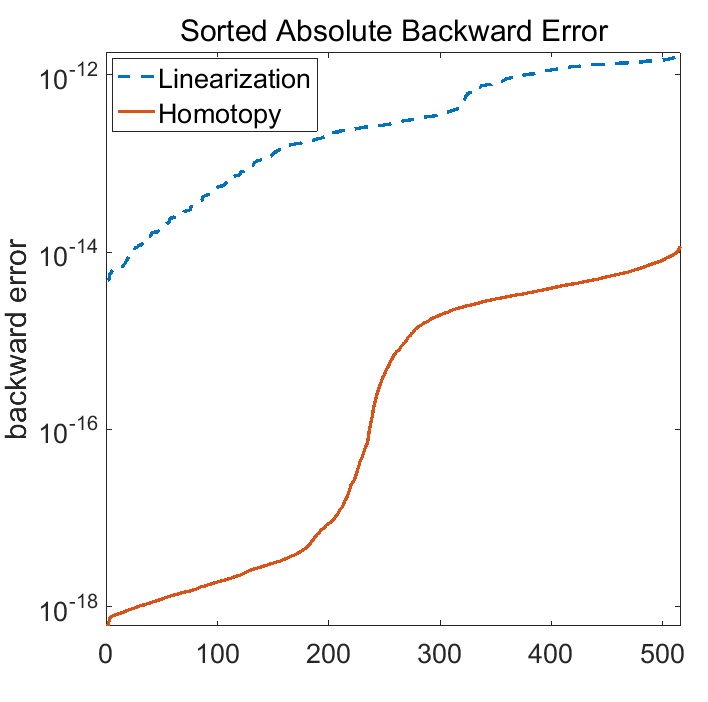}
\includegraphics[scale=0.43]{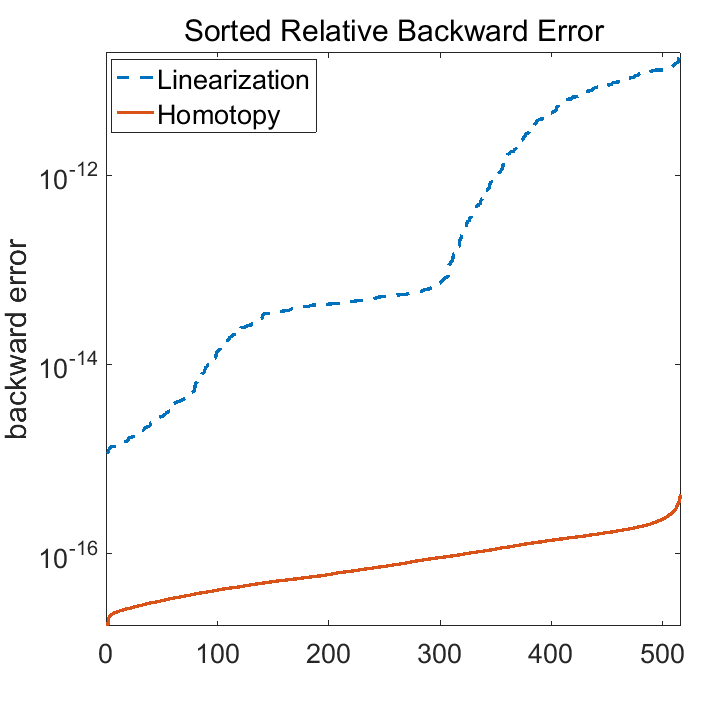}
\label{ex2_BErr}
\end{figure}

The dimension of this PEP is fixed at $n=129$. 
As before, we compare the accuracy and timings of homotopy and linearization methods on this problem. We present the absolute and relative backward errors for all computed eigenpairs in \cref{ex2_BErr} and tabulate the best, average, median, worst performance in \cref{Ave_berr_APP2}.  The elapsed timings are given in  \cref{timing_APP2}. All speed and accuracy tests are run ten times with  the best, average, median and worst results recorded.


\begin{table}[H]
\footnotesize
\centering
\caption{\textsc{Planar waveguide problem --- absolute and relative backward errors.} Best, average, median, worst absolute  (top table) and relative (bottom table) backward errors of all eigenpairs computed via homotopy and linearization methods.}
\label{Ave_berr_APP2}
\begin{tabular}{ c| c| c}
\hline
 $n=129$   &LINEARIZATION ABS.\ BK.\ ERR.& HOMOTOPY ABS.\ BK.\ ERR.\\
\hline
BEST & 4.77790E-15 & 6.06981E-19 \\
MEAN & 5.38568E-13 & 2.05656E-15	 \\
MEDIAN & 2.83418E-13 & 6.90006E-16	 \\
WORST & 1.76856E-12 & 1.166201E-14	 \\
\hline
$n=129$   &LINEARIZATION REL.\ BK.\ ERR.& HOMOTOPY REL BK.\ ERR.\\
\hline
BEST & 1.14289E-15 & 1.72613E-17 \\
MEAN & 2.53466E-12 & 9.48695E-17	 \\
MEDIAN & 5.37243E-14 & 7.64027E-17	 \\
WORST & 1.87416E-11 & 4.24275E-16
\end{tabular}
\end{table}

\begin{table}[H]
\footnotesize
\centering
\caption{\textsc{Planar waveguide problem --- speed.} 
Elapsed timings (in seconds) for homotopy and linearization methods. The homotopy method is run in parallel on $80$ cores.}
\label{timing_APP2}
\begin{tabular}{ c| c| c}
$n=129$   &LINEARIZATION TIMINGS& HOMOTOPY TIMINGS\\
\hline
BEST & 2.6670 & 2673.0350 \\
MEAN & 2.7321 & 3013.9718	 \\
MEDIAN & 2.7213 & 2973.6600\\
WORST & 2.8414 & 3460.5290
\end{tabular}
\end{table}

The results obtained for this quartic PEP arising from the planar waveguide problem is consistent with what we have observed for the acoustic wave QEP in \cref{ss:QEPEx} as well as the  randomly generated PEPs and QEPs in \cref{Numerical results} --- while homotopy method requires much longer running times  than the linearization method, its results are  also vastly superior in terms of accuracy.

In summary, if our main goal is to obtain accurate solutions to polynomial eigenvalue problems, particularly when all eigenpairs  are needed, then expending additional resources (more cores and longer computing time) to employ the homotopy method is not only worthwhile but perhaps inevitable --- we know of no other alternative that would achieve the same level of accuracy.

\section*{Acknowledgment}

The work in this article is generously supported by DARPA D15AP00109 and NSF IIS 1546413. LHL is supported by a DARPA Director's Fellowship.  JIR is supported by 
a University of Chicago Provost Postdoctoral Scholarship.

\bibliography{BIB_AccuratePEP}
\bibliographystyle{siamplain}

\end{document}